\documentclass [12 pt, letterpaper] {article}

        \usepackage {times,amssymb,amsmath,array,calc,enumerate}
   {\begin {list}{} 
   {\usecounter{example}%
   \small%
   \setlength{\labelsep}{0pt}\setlength{\leftmargin}{40pt}%
   \setlength{\labelwidth}{0pt}%
   \setlength{\listparindent}{0pt}}}%
   {\end{list}}

\renewcommand {\theequation} {\@arabic\c@equation}

\makeatother
   {\begin {list}{} 
   {
   \small%
    \setlength{\itemindent}{10pt}
    \setlength{\labelsep}{5pt}\setlength{\leftmargin}{40pt}%
    \setlength{\labelwidth}{5pt}%
   \setlength{\itemsep}{-1.8pt}
    \setlength{\listparindent}{0pt}
   }
   }%
   {\end{list}}

\newcommand{\pb}{\vspace{0.3cm} \noindent}
\newcommand{\N}{\mathcal {N}}
\newcommand{\A}{\mathcal{A}}
\newcommand{\M}{\mathcal{M}}

\newcommand{\mb}{\mathbf}

\renewcommand {\theequation} {\arabic{equation}}

\newtheorem{Def}{Definition}[section]
\newtheorem{Th}{Theorem}[section]
\newtheorem{Col}{Corollary}[section]
\newtheorem{Prop}{Proposition}[section]
\newtheorem{lemma}{Lemma}[section]
\begin {document}

\centerline{\Large On Voiculescu's Semicircular Matrices }

\thispagestyle{empty} \vspace{2cm}

\centerline{\large Liming Ge \qquad and \qquad  Junhao Shen}

\vspace{1cm}

\centerline{Mathematics Department, University of New Hampshire,
Durham, NH, 03824}

\vspace{0.5cm}

\centerline{email: liming@math.unh.edu \quad and \quad
jog2@cisunix.unh.edu}

\vspace{3cm}

\pb{\bf Abstract:} Assume $\N$ is a von Neumann algebra of type
II$_1$ with a tracial state $\tau_{\N}$, and $\M$ is the von Neumann
algebra of the $n\times n$ matrices over $\N$ with the canonical
tracial state $\tau_{\M}$. Let $\mathcal D_n$ be the subalgebra of
$\M$ consisting of   scalar diagonal matrices in $\M$. In this
article, we study the properties of   semicircular elements in $\M$
that are free from $\mathcal D_n$ with respect to $\tau_{\M}$. Then
we define a new concept ``matricial distance" of two elements in
$\M$ and compute the matricial distance between two free
semicircular elements in $\M$.

\newpage

\section{Introduction}

The theory of free probability was developed by Voiculescu in early
80's in the last century. In his influential paper [11], Voiculescu
introduced   concepts of semicircular and circular systems (in this
paper  these systems are called as    ``semicircular and circular
matrices") and used them to obtain the surprising relationship
between free group factors. One of essential arguments in [11] is
that the semicircular matrices can be obtained from the matrices
over another free probability space if the entries of the matrices
are, in an appropriate way, chosen to be semicircular or circular
and free. Moreover, the semicircular matrices obtained in this way
are free with the algebra generated by diagonal matrices over
scalars (see Theorem 2.1).

More specifically, let $(\N, \tau_{\N})$ be a free probability space
and $\M=\N\otimes M_n(\Bbb C)$ be the algebra of $n\times n$
matrices over $\N$ with the canonical tracial state $\tau_{\M}$. Let
$\{e_{ij}\}_{1\le i,j\le n}$ be the canonical system of matrix units
of $M_n(\Bbb C)$. Let $\mathcal D_n$ be the subalgebra of $\M$
generated by $\{I_{\N}\otimes e_{ii}\}_{1\le i\le n}$. Let
$A=\sum_{1\le i,j\le n} a_{ij}\otimes e_{ij}$ be a semicircular
element of $\M$ such that $\{a_{ij}\}_{1\le i,j\le n}$ are, in an
appropriate way, chosen to be semicircular or circular and free in
$\N$. (Such carefully constructed semicircular matrix
 will be called as    Voiculescu's semicircular matrix.) (see Def. 2.5) Then
Voiculescu proved his celebrated result in [11] that $A$ and
$\mathcal D_n$ are free with respect to $\tau_{\M}$.

The concept of semicircular element (or matrix) has now been a key
concept in the theory of free probability and its application on von
Neumann algebras.   There are many further studies on semicircular
elements after Voiculescu's paper. For example, the concept of
R-cyclic matrices, as a generalization of Voiculescu's semicircular
and circular matrices, was introduced
 in a remarkable paper
[6]. It was proved there that  a matrix $A$, $n\times n$ matrix over
$(\N, \tau_{\N})$,  is R-cyclic if and only if $A$   is free from
$M_n(\Bbb C)$ with amalgamation over $\mathcal D_n$. Some examples
of R-cyclic matrices are given there.

The purpose of the paper is to further discuss some special
properties of a semicircular element $B$ in $\M$  when $B$ and
$\mathcal D_n$ are free in $\M$ with respect to $\tau_{\M}$. The
first result obtained  is Theorem 3.3, an inverse statement of
Vouculecu's celebrated result (see also Theorem 2.1) in [11].

\vspace{0.3cm} \noindent {\bf Theorem:}
  Suppose $B= \sum_{1\le i,j\le n} b_{ij}\otimes e_{ij}$ is a
  semicircular element of $(0,1)$ in $\mathcal M$. If $B$ and $\mathcal
  D_n$ are free with respect to $\tau_{\mathcal M}$, then there is a family of unitary elements
  $\{u_1, \ldots, u_n\}$ in $\mathcal N$ such that $U^*BU$ is a
  Voiculescu's semicircular matrix with $U=\sum_{1\le
  i\le n} u_i\otimes e_{ii}$.
\vspace{0.3cm}

 Using
the preceding theorem, we are able to give   characterization of
semicircular elements that are free with $\mathcal D_n$ in $\M$. In
more details, we have

 \vspace{0.3cm}
\noindent {\bf Theorem: }Suppose $B= \sum_{1\le i,j\le n}
b_{ij}\otimes e_{ij}$ is a
  self-adjoint element in $\mathcal M$ such that   $B$ and $\mathcal
  D_n$ are free with respect to $\tau_{\M}$. Then the following are equivalent.
  \begin{enumerate}[(i)]
      \item $B$ is a semicircular element in $\mathcal M$.
      \item There is a family of unitary elements
  $\{u_1, \ldots, u_n\}$ in $\mathcal N$ such that $U^*BU$ is a
  Voiculescu's   semicircular matrix, where $U=\sum_{1\le
  i\le n} u_i\otimes e_{ii}$.
  \item For all $i , j , \ 1\le i  \ne j \le
n$, $b_{ij}^*b_{ij} $ and  $b_{jj}$ are free   with respect to
$\tau_{\mathcal N}$.
\item There are some  $i_0, j_0, \ 1\le i_0 \ne j_0\le
n$, $b_{i_0j_0}^*b_{i_0j_0} $ and  $b_{j_0j_0}$ are free    with
respect to $\tau_{\mathcal N}$.
\item For all $i , j , \ 1\le i  \ne j \le
n$, $b_{ij}b_{ij}^*$ and  $b_{ii}$ are free    with respect to
$\tau_{\mathcal N}$.
\item There are some  $i_0, j_0, \ 1\le i_0 \ne j_0\le
n$, $b_{ i_0j_0} b_{i_0j_0}^*$ and  $b_{i_0i_0}$ are free    with
respect to $\tau_{\mathcal N}$.
  \end{enumerate}
Therefore, under the conjugation of  diagonal unitary matrices,
Voiculescu's semicircular matrix is the unique form of semicircular
elements that are free from $\mathcal D_n$  with respect to
$\tau_{\M}$. Moreover, under the assumption that the self-adjoint
element $B$ is free from $\mathcal D_n$ in $\M$, a very weak
condition on   freeness among entries of $B$ (for example an off
diagonal entry $a_{i_0j_0}$ and a diagonal entry $a_{j_0j_0}$ in
same column are free with respect to $\tau_{\N}$), will imply that
$B$ is a semicircular element (also see Theorem 3.1 and 3.2). In
this sense, Voiculescu's semicircular matrix is the only
self-adjoint matrix that  is free from $\mathcal D_n$ in $\M$ and
has free entries inside.

 The second half of this paper is  devoted to   study    another aspect of
 semicircular elements, which is inspired by [9], a remarkable paper by S. Popa.
In [9], by  showing that any standard generator $u_g$  of free group
factors $L(F(n))$ is not contained in any hyperfinite II$_1$
subfactor of $L(F(n))$ (also see [2]), S. Popa answered an open
question asked by R. Kadison (see [4]). Since the hyperfinite II$_1$
factor can be ``approximated" by large sized type I$_k$ factors, S.
Popa's result can be interpreted in the following way. Suppose that
$u_g$ and $u_h$ are two different standard generators of $L(F(n))$
(hence they are free with each other). Once $u_g$ can be almost
``contained" in a ``large" type I$_k$ factor $M_k$, then we would
expect that $u_h$ will be far away from $M_k$, i.e. the distance
between $u_h$ and $M_k$ is going to be large. We define a new
concept called by ``matricial distance" to describe this asymptotic
phenomenon of the distance between $u_h$ and $M_k$ that almost
contains $u_g$. More generally, assume that $A$, $B$ are two free
semicircular elements of $(0,1)$  in a type II$_1$ factor $\M$.
(Note, here, we do not require   $\M$ to be a free group factor).
 For every   positive integer $k$ and every family of mutually orthogonal
 equivalent projections $P_1,\ldots, P_k$  with sum $I$ in $\mathcal
 A$, the abelian von Neumann subalgebra generated by $A$ in $\M$, let $\mathcal M_k$ be any type I$_k$ factor in $\mathcal M$  such that the
 diagonal projections of $M_k$ are $\{P_i\}_{1\le i\le k}$.  Let $E_k $ be the conditional
 expectation from $\M$ onto $M_k$. Then we showed in Proposition 4.1 that $||E_k(B)||_2\le 7/\sqrt[8]
 {k},$ or   the
matricial distance between two free semicircular elements in   $\M$
is equal to 1.

Further calculation of   matricial distance between   two free
self-adjoint elements, not required to be semicircular, will be
carried out in our forthcoming paper.

 The organization of the paper is as follows. The basic
 knowledge is reviewed in section 2. Section 3 is devoted to
 prove   ``uniqueness" of Voiculescu's semicircular matrix. The necessary  and  sufficient  conditions for  a matrix to be semicircular and
 free with $\mathcal D_n$ are given in Theorem 3.5. As a
 corollary, the result of ``free compression" in [7]   is
 reproved in the same section.     Proposition 4.1  is proved in section 4.

\section{Preliminary}

 A  {\em  free probability space} is a pair
    $(\mathcal M, \tau)$ where $\mathcal M$ is a von Neumann algebra
    and $\tau_{\mathcal M}$ is a normal state ([5]). We assume that $\mathcal M$ has a
    separable predual and $\tau$ is a faithful normal tracial state ([5]).
    (So that $\mathcal M$ is a finite von Neumann algebra with a trace). Elements
    of $\mathcal M$ are called non-commuting {\em random variables}.

    \begin{Def}{  ([14]) The {\em distribution} of a
    random variable $A$ in $(\mathcal M, \tau)$ is a linear functional
    $\mu$ on $\Bbb C[x]$, the polynomial ring with variable $x$ and
    coefficients in $\Bbb C$, such that $\mu(\psi(x))=\tau (\psi(A))$
    for all $\psi(x)$ in $\Bbb C[x]$.}
    \end{Def}

    \begin{Def} {A {\em semicircular element  }  $A$  in $\mathcal M$ is one whose
    distribution $\gamma_{a,r}$ satisfies the
     {\em ``semicircle law"}, one whose distribution
     $\gamma_{a,r}$ (centered at $a$ in $\Bbb R$ with radius
     $r>0$) maps $\Bbb C[x]$ to $\Bbb C$ according to the
     equation:
    $$
    \gamma_{a, r} (\psi)=\frac 2 {\pi r^2}\int_{a-r}^{a+r} \psi(t)
    \sqrt{r^2-(t-a)^2} dt.
    $$
     Such $A$ is called a {\em  semicircular element of $(a,r)$.}

     \vspace{0.3cm} The positive element $H$ in $\mathcal M$ is called a {\em quarter-circular  element of $(a,r)$ } if   there
     is a       semicircular element of $(a,r)$, $B$, in $\M$   such that $H= \sqrt{B^*B}$.
      }
    \end{Def}
    \begin{Def}{
    ([14]) The {\em joint distribution} of a family of
    random variables $A_i$, $i \in \mathcal I$, in $(\mathcal M,
    \tau)$ is a linear functional $\mu$ on
       $\Bbb C \langle x_i, i \in
    \mathcal I\rangle$, the noncommutative polynomial ring with noncommuting
    variables $x_i$, such that $\mu(\psi(x_{i_1}, \ldots,
    x_{i_n}))=\tau (\psi(A_{i_1}, \ldots, A_{i_n}))$ for every $\psi$
    in $\Bbb C \langle x_i, i \in \mathcal I \rangle $.}
    \end{Def}

    \begin{Def}{
    ([14]) The von Neumann subalgebras $\mathcal M_i,
    i \in \mathcal I$ of $\mathcal M$ are {\em free} with respect to
    the trace $\tau$ if $\tau (A_1\ldots A_n)=0$ whenever $A_j \in
    \mathcal M_{i_j}, i_1\ne \ldots \ne i_n$ and $\tau (A_j)=0$ for
    $1 \le j \le n$ and every $n$ in $\Bbb N$.
     (Note that $i_1$ and $i_3$, for example, may be equal:
     ``adjacent" $A_i$s are not in the same $\mathcal M_i$). A
     family of self-adjoint elements $\{A_1, \ldots, A_n\}$ is
     free with respect to the trace $\tau$ if the von Neumann
     subalgebras $\mathcal M_i$ generated by the $A_i$ are
     free with respect to the trace $\tau$.

    }
    \end{Def}

Let $\N$ be a von Neumann algebra with a tracial state $\tau_{\N}$,
and $\M$ be  $\mathcal N\otimes M_n(\Bbb C)$ for   some integer $n$
in $\Bbb N$.
  Let $\{e_{ij}\}_{1\le i,j\le n}$ be the canonical system of matrix
  units of $M_n(\Bbb C)$ in $\mathcal M$. Let $\mathcal D_n$ be the
  abelian von Neumann subalgebra generated by $\{I_{\N}\otimes e_{ii}\}_{1\le i\le n}$ in $\mathcal M$.
  Let $\tau_{\M}$ be the canonical trace on $\M$ defined as:
  $$\tau_{\M}(A) = \frac 1 n \sum_{i=1}^n\tau_{\N}(a_{ii}),$$
  for every $A =\sum_{1\le i,j\le n} a_{ij}\otimes e_{ij}$ in $\mathcal
  M$.

\begin{Def} A self-adjoint element $A =\sum_{1\le i,j\le n} a_{ij}\otimes e_{ij}$ in $\mathcal M$ is called {\em Voiculecu's semicircular
matrix} if the following hold:
\begin{enumerate} [(a)]
\item $\{Re\  a_{
i,j }, Im\  a_{ i,j } \  |\ 1 \le i\le j \le n\} $ is a family of
free elements in $\mathcal N$ with respect to $\tau_{\mathcal N}$
    \item   $ Re\  a_{
i,j }   $ and $ Im\  a_{ i,j }   $  are   semicircular elements
    of
    $(0, \frac 1 { \sqrt{2n}})$ in $\mathcal N$   with respect to
     $\tau_{\mathcal N}$,
    for $1 \le i<j \le
n-1$.
  \item  Each $    a_{
i,n }   $ is a quarter-circular element of
    $(0, \frac 1 {\sqrt {n}})$ in $\mathcal N$ with respect to
      $\tau_{\mathcal N}$,
    for $1 \le i<
n-1$.
   \item Each $ a_{
j,j } $ is a semicircular element of  $(0, \frac 1 { \sqrt{n}})$
   in $\mathcal N$ with respect to
    $\tau_{\mathcal N}$,
   for $1 \le j\le n$.
\end{enumerate}
\end{Def}

\begin{Def} A  family of self-adjoint elements $\{A_{k} =\sum_{1\le i,j\le n} a_{ij}^{(k)}\otimes e_{ij}\}_{1\le k\le m}$ in $\mathcal M$ is called
{\em a standard family of Voiculecu's  semicircular matrices} if the
following hold:
\begin{enumerate}[(a)]
\item $\{Re\  a_{
i,j }^{(k)}, Im\  a_{ i,j }^{(k)} \  |\ 1 \le i\le j \le n, 1\le
k\le m\} $ is a family of free elements in $\mathcal N$ with respect
to $\tau_{\mathcal N}$
    \item   $ Re\  a_{
i,j }^{(k)}  $ and $ Im\  a_{ i,j }^{(k)}  $   are   semicircular
elements
    of
    $(0, \frac 1 { \sqrt{2n}})$ in $\mathcal N$   with respect to
     $\tau_{\mathcal N}$,
    for $1 \le i<j \le
n-1, 1\le k\le m$.
  \item  Each $     a_{
i,n }^{(1)}     $ is a quarter-circular element of
    $(0, \frac 1 {\sqrt {n}})$ in $\mathcal N$ with respect to
      $\tau_{\mathcal N}$,
    for $1 \le i<
n-1$.
 \item    $ Re\  a_{
i,n }^{(k)}  $ and $ Im\  a_{ i,n }^{(k)}  $  are semicircular
elements
    of
    $(0, \frac 1 { \sqrt{2n}})$ in $\mathcal N$   with respect to
     $\tau_{\mathcal N}$,
    for $1 \le i  \le
n-1, 2\le k\le m$.
   \item Each $ a_{
i,i }^{(k)} $ is a semicircular element of  $(0, \frac 1 {
\sqrt{n}})$
   in $\mathcal N$ with respect to
    $\tau_{\mathcal N}$,
   for $1 \le i\le n, 1\le k\le m$.
\end{enumerate}
\end{Def}

\vspace{0.3cm}\noindent  Following the notations as above,
Voiculescu proved the following remarkable theorem.
\begin{Th} (From [11]) Suppose $\{A_{k}\}_{1\le k\le m}$ is { a standard family
of Voiculescu's  semicircular matrices}. Then $\{\mathcal D_n,
A_{1}, \ldots, A_{m}\}$ are free with respect to $\tau_{\M}$.

\end{Th}

\vspace{2cm}

\section{Uniqueness of Voiculescu's Semicircular Matrix}
Following the notations from preceding section, we let $\N$ be a von
Neumann algebra with a tracial state $\tau_{\N}$, and $\M$ be
$\mathcal N\otimes M_n(\Bbb C)$ for   some integer $n$ in $\Bbb N$
with the canonical tracial state $\tau_{\M}$. Let
$\{e_{ij}\}_{i,j=1}^{ n}$ be the canonical system of matrix units of
$M_n(\Bbb C)$ in $\mathcal M$.

  Our next proposition shows that if two self-adjoint elements
share certain properties of freeness and   agree on first two
moments then they have the same distributions.
 \begin{Prop}
Suppose that $B=\sum_{1\le i,j\le n} b_{ij}\otimes e_{ij}$,
$X=\sum_{1\le i,j\le n} x_{ij}\otimes e_{ij}$ are two self-adjoint
elements in $\mathcal M$ (with $b_{ij}$'s and $x_{ij}$'s in
$\mathcal N$) such that, for some $i_0, j_0, \ 1\le i_0 \ne j_0\le
n$, the following hold.
\begin{enumerate}[(a)]
\item $B$ and $\{I_{\N}\otimes e_{i_0i_0}, I_{\N}\otimes  e_{j_0j_0}\}$  are free   with respect to
 $\tau_{\mathcal M}$;
\item $X$ and $\{I_{\N}\otimes e_{i_0i_0},I_{\N}\otimes  e_{j_0j_0}\}$  are free  with respect to
 $\tau_{\mathcal M}$;
\item   $b_{j_0i_0}b_{i_0j_0}$ $(= b_{i_0j_0}^*b_{i_0j_0} )$ and  $b_{j_0j_0}$ are free  with respect to
 $\tau_{\mathcal N}$;
\item  $x_{j_0i_0}x_{i_0j_0}$ $(= x_{i_0j_0}^*x_{i_0j_0})$ and  $x_{j_0j_0}$ are free  with respect to
 $\tau_{\mathcal N}$;
 \item  $\tau_{\mathcal
M}(B)= \tau_{\mathcal M}(X)$ and $\tau_{\mathcal M}(B^2)=
\tau_{\mathcal M}(X^2)$.
\end{enumerate}
Then $\tau_{\mathcal M}(B^m)= \tau_{\mathcal M}(X^m)$ for every
$m\ge 1$.
 \end{Prop}

\noindent{\bf Proof:} We need only to prove the proposition under
the assumption that $ \tau_{\mathcal M}(B)=0,\ \tau_{\mathcal
M}(B^2)=1.$ Otherwise we let $\tilde B= \frac 1 {\|B-\tau_{\mathcal
M}(B)\|_2} (B-\tau_{\mathcal M}(B))$ and $\tilde X= \frac 1
{\|X-\tau_{\mathcal M}(X)\|_2} (X-\tau_{\mathcal M}(X))$. Then
$\tilde B$ and $\tilde X$ are self-adjoint elements that satisfy all
conditions in the proposition and $ \tau_{\mathcal M}(\tilde
B)=\tau_{\mathcal M}(\tilde X)=0,\ \tau_{\mathcal M}(\tilde
B^2)=\tau_{\mathcal M}(\tilde X^2)=1.$ The result that
$\tau_{\mathcal M}(\tilde B^m) = \tau_{\mathcal M}(\tilde X^m) $ for
every $m\ge 1$ will imply that $\tau_{\mathcal M}( B^m) =
\tau_{\mathcal M}( X^m) $.

 \vspace{0.3cm}  If it brings no confusion, we will write $I_{\N}\otimes
e_{ij}$ as $e_{ij}$.

\vspace{0.3cm}

\noindent{\bf Claim I:} $\tau_{\mathcal
  N}(b_{ij}b_{ji}) = \tau_{\mathcal
  N}(x_{ij}x_{ji}), \qquad \forall i,j \in
\{i_0, j_0\}.$

\vspace{0.3cm}
  Note that
 we have, for any $i\in
\{i_0, j_0\}, j \in \{i_0, j_0\}$  and any $A$ in $\mathcal M$,
$$  \begin{aligned}\tau_{\mathcal
  M}(e_{i i }  Ae_{j j } A)  
 & =\tau_{\mathcal
  M}((e_{i i }-\frac 1 n  )  A(e_{j j }-\frac 1 n  )
  A)  -\frac 1 n \tau_{\mathcal
  M}(  (e_{j j }-\frac 1 n )
  AA)\\
  &\quad   -\frac 1 n \tau_{\mathcal
  M}(  (e_{i i }-\frac 1 n )
  AA)   +\frac 1 {n^2} \tau_{\mathcal
  M}( AA )
\end{aligned}
  $$
It follows from that facts $\tau_{\mathcal M}(B)= \tau_{\mathcal
M}(X)=0$, $\tau_{\mathcal M}(B^2)= \tau_{\mathcal M}(X^2)=1$ and
conditions (a), (b), that
\begin{equation}
  \tau_{\mathcal
  M}(e_{i i }  Be_{j j } B) = \tau_{\mathcal
  M}(e_{i i }  Xe_{j j } X) \qquad \forall i,j \in
\{i_0, j_0\}.\end{equation} i.e.
$$
  \tau_{\mathcal
  N}(b_{ij}b_{ji}) = \tau_{\mathcal
  N}(x_{ij}x_{ji}) \qquad \forall i,j \in
\{i_0, j_0\}.$$
 \vspace{0.3cm}

Instead of  proving $\tau_{\mathcal M}(B^m)= \tau_{\mathcal M}(X^m)$
for every $m\ge 1$, we are going to prove a stronger result.

\vspace{0.3cm}

 \noindent {\bf Claim II:} For
each positive integer $m$, we have
\begin{enumerate}[(i)]
    \item $\tau_{\mathcal
  M}(e_{i_0i_0} (Be_{j_0j_0})^{m-1}B)
  =\tau_{\mathcal
  M} (e_{i_0i_0} (Xe_{j_0j_0})^{m-1}X) $.
  \item $\tau_{\mathcal N} ((b_{j_0j_0} )^m)
  =\tau_{\mathcal N}((x_{j_0j_0} )^m)$.
  \item   For $2\le l\le m,$ $1 \le t_1, \ldots, t_l \le m$ with $  t_1+\cdots + t_l\le m+1,$   and
  $s_1,\ldots, s_l \in \{i_0, j_0\}$, we have \quad
   $\tau_{\mathcal M}((e_{s_1s_1}-\frac 1 n)B^{t_1}\ldots
  (e_{s_ls_l}-\frac 1 n)B^{t_l})= \tau_{\mathcal
  M}((e_{s_1s_1}-\frac 1 n)X^{t_1}\ldots
  (e_{s_ls_l}-\frac 1 n)X^{t_l}). $
  \item $\tau_{\mathcal M}(e_{i_0i_0}(Be_{j_0j_0})^mB)=\tau_{\mathcal
  M} (e_{i_0i_0}(Xe_{j_0j_0})^mX) $
  \item $\tau_{\mathcal M}(B^{m+1})=\tau_{\mathcal
  M} (X^{m+1})$\end{enumerate}

\vspace{0.3cm} The claim will be proved by the induction on $m$.

\vspace{0.3cm}

\noindent {\bf The case when $\mb m=\mb 1$:} $(i), (ii)$ are
directly from conditions (a) and (b). $(iv)$ is  equation (1). $(v)$
is trivial.

\vspace{0.3cm}

\noindent {\bf The case when $\mb m=\mb 2$:} $(i), (ii)$ are from
equation (1) directly. For any $1\le t_1, t_2 \le 2$, $t_1+t_2\le 3$
and $\{s_1,s_2\}\subset \{i_0, j_0
  \}$, we have  that$$  \begin{aligned}&\tau_{\mathcal
  M}((e_{s_1 s_1 } -\frac 1 n  ) B^{t_1}(e_{s_2s_2 }-\frac 1 n  ) B^{t_2}) \\&  = \tau_{\mathcal
  M}((e_{s_1 s_1  }-\frac 1 n )(  B^{t_1}-\tau_{\mathcal
  M} (B^{t_1})+\tau_{\mathcal
  M} (B^{t_1}))(e_{s_2s_2 }-\frac 1 n  )
  (  B^{t_2}-\tau_{\mathcal
  M} (B^{t_2})+\tau_{\mathcal
  M} (B^{t_2}))) \\
 & =\tau_{\mathcal
  M}((e_{s_1 s_1  }-\frac 1 n )(  B^{t_1}-\tau_{\mathcal
  M} (B^{t_1}) )(e_{s_2s_2 }-\frac 1 n  )
  (  B^{t_2}-\tau_{\mathcal
  M} (B^{t_2}) )) \\
&\quad + \tau_{\mathcal
  M}((e_{s_1 s_1  }-\frac 1 n )(  \tau_{\mathcal
  M} (B^{t_1}))(e_{s_2s_2 }-\frac 1 n  )
  (  B^{t_2}-\tau_{\mathcal
  M} (B^{t_2}) )) \\
  &\quad +\tau_{\mathcal
  M}((e_{s_1 s_1  }-\frac 1 n )(  B^{t_1}-\tau_{\mathcal
  M} (B^{t_1}) )(e_{s_2s_2 }-\frac 1 n  )
  (  \tau_{\mathcal
  M} (B^{t_2}))) \\
&\quad +\tau_{\mathcal
  M}((e_{s_1 s_1  }-\frac 1 n )(  \tau_{\mathcal
  M} (B^{t_1}))(e_{s_2s_2 }-\frac 1 n  )
  (  \tau_{\mathcal
  M} (B^{t_2}) )) \\
  &  =    
  \tau_{\mathcal
   M} (B^{t_1})  \tau_{\mathcal
   M} (B^{t_2}) \tau_{\mathcal
  M}((e_{s_1 s_1  }-\frac 1 n )(e_{s_2s_2 }-\frac 1 n  ))\qquad  \qquad
  \qquad (\text {because of condition (a)})
    \end{aligned} $$
   Similar argument shows that
   $$ \begin{aligned}
&\tau_{\mathcal
  M}((e_{s_1 s_1 } -\frac 1 n  ) X^{t_1}(e_{s_2s_2 }-\frac 1 n  )
  X^{t_2}) =
  \tau_{\mathcal
  M}(X^{t_1})   \tau_{\mathcal
  M} (X^{t_2})  \tau_{\mathcal
  M}((e_{s_1 s_1  }-\frac 1 n )(e_{s_2s_2 }-\frac 1 n  )
   )
\end{aligned}$$
Note that  $1\le t_1,t_2\le 2$ and $\tau_{\M}(B)=\tau_{\M}(X)$,
$\tau_{\M}(B^2)=\tau_{\M}(X^2)$. We have that
   $(iii)$ holds for $m=2$.

  As for $(iv)$, we know   that
 $$\begin{aligned}
    \tau_{\mathcal M}(e_{i_0i_0}(Be_{j_0j_0})^2B)&=
    \frac 1 n \tau_{\mathcal N} (b_{i_0j_0}  b_{j_0j_0}
    b_{j_0i_0})\\& = \frac 1 n \tau_{\mathcal N}   (b_{j_0j_0}
     ) \tau_{\mathcal N} (
    b_{j_0i_0}b_{i_0j_0})\qquad \quad  (\text {because of condition (c)})
   \\
&= \frac 1 n \tau_{\mathcal N}  (x_{j_0j_0}
     ) \tau_{\mathcal N}(
    x_{j_0i_0}x_{i_0j_0}) \qquad \quad  (\text {because of equation (1)})\\
    &   =
    \frac 1 n \tau_{\mathcal N} (x_{i_0j_0}  x_{j_0j_0}
    x_{j_0i_0})    \qquad \quad \qquad  (\text {because of condition (d)})\\
   &=\tau_{\mathcal M}(e_{i_0i_0}(Xe_{j_0j_0})^2X)
 \end{aligned}$$

 Now we are ready to prove $(v)$. We have
 $$
 \begin{aligned}
   \tau_{\mathcal M}&(e_{i_0i_0}(Be_{j_0j_0})^2B)\\ &= \tau_{\mathcal M}((e_{i_0i_0}-\frac 1 n +\frac 1 n)B(e_{j_0j_0}-\frac 1 n
  +\frac 1 n)B(e_{j_0j_0}-\frac 1 n +\frac 1 n)B)\\
&= \tau_{\mathcal M}((e_{i_0i_0}-\frac 1 n  )B(e_{j_0j_0}-\frac 1 n
   )B(e_{j_0j_0}-\frac 1 n  )B)\\
  &\quad + \frac 1 n  \tau_{\mathcal M}(  B(e_{j_0j_0}-\frac 1 n
  )B(e_{j_0j_0}-\frac 1 n  )B)
              + \frac 1 n \tau_{\mathcal M}((e_{i_0i_0}-\frac 1 n  )B
    B(e_{j_0j_0}-\frac 1 n  )B)\\
&\quad +  \frac 1 n\tau_{\mathcal M}((e_{i_0i_0}-\frac 1 n
)B(e_{j_0j_0}-\frac 1 n
  )BB) +   \frac 1 {n^2} \tau_{\mathcal M}(   B
    B(e_{j_0j_0}-\frac 1 n )B)\\
&\quad + \frac 1 {n^2}\tau_{\mathcal M}((e_{i_0i_0}-\frac 1 n  )B
 B   B) + \frac 1 {n^2} \tau_{\mathcal M}(
B(e_{j_0j_0}-\frac 1 n
  )B   B)\\ &\quad    +\frac 1 {n^3} \tau_{\mathcal M}(   B    B
  B)\\
  &=\frac 1 {n^3}\tau_{\mathcal M}( B^3) \qquad \qquad\qquad \qquad \qquad \qquad \qquad \qquad \quad  (\text {because of condition (c)})
 \end{aligned}
 $$
On the other hand, similar computation shows that
$$
    \tau_{\mathcal M}(e_{i_0i_0}(Xe_{j_0j_0})^2X) = \frac 1 {n^3} \tau_{\mathcal M} (
    X^3)
 $$
Combining with $(iv)$, we know that $\tau_{\mathcal M}(B^3)
=\tau_{\mathcal M}(X^3) $.

\vspace{0.5cm}

\noindent {\bf The case when $\mb m=\mb k+\mb 1$:}  Assume that
$(i)$, to $(v)$ hold when $m=k$. Consider
  $m=k+1$. We have
  $$
\begin{aligned}
  \tau_{\mathcal
  M}(e_{i_0i_0} (Be_{j_0j_0})^{k}B) &= \frac 1 n \tau_{\mathcal N}
  (b_{i_0j_0}(b_{j_0j_0})^{k-1}b_{j_0i_0})\\
  &=\frac 1 n \tau_{\mathcal N}(b_{j_0i_0}b_{i_0j_0}) \tau_{\mathcal
  N}((b_{j_0j_0})^{k-1}) \\
  &= \frac 1 n\tau_{\mathcal N}(x_{j_0i_0}x_{i_0j_0}) \tau_{\mathcal
  N}((x_{j_0j_0})^{k-1})\\
 & \qquad \qquad \qquad    \ \ \text {(because of $(ii)$ of induction hypothesis) }
  \\
&= \frac 1 n \tau_{\mathcal N}
   (x_{i_0j_0}(x_{j_0j_0})^{k-1}x_{j_0i_0})\\
  &= \tau_{\mathcal M}(e_{i_0i_0} (Xe_{j_0j_0})^{k}X)
\end{aligned}
  $$
So, $(i)$ holds when $m=k+1$. From $(v)$ of induction hypothesis, we
know that $\tau_{\mathcal M}(B^i)=\tau_{\mathcal M} (X^i)$ for $1\le
i\le k+1$. Since $B$ and $e_{i_0i_0}$ are free, we have
$$ \tau_{\mathcal M} ((e_{i_0i_0}B)^{k+1})=\tau_{\mathcal M}((e_{i_0i_0}X)^{k+1}). $$
It follows  that $\tau_{\mathcal N}((b_{i_0i_0})^{k+1})
=\tau_{\mathcal M}((x_{i_0i_0})^{k+1})$. Therefore $(ii)$ holds.

\vspace{0.3cm}

\noindent As for $(iii)$,  when $2\le l\le k+1,$ $1 \le t_1, \ldots,
t_l \le k+1$ with $ t_1+\cdots + t_l\le k+2,$   and
  $s_1,\ldots, s_l \in \{i_0, j_0\}$, we have
  $$\begin{aligned}
  &\tau_{\mathcal M}((e_{s_1s_1}-\frac 1 n)B^{t_1}\ldots
  (e_{s_ls_l}-\frac 1 n)B^{t_l})\\
&= \tau_{\mathcal M}((e_{s_1s_1}-\frac 1 n)(B^{t_1}-\tau_{\mathcal
M}(B^{t_1})+ \tau_{\mathcal M}(B^{t_1}))\ldots
  (e_{s_ls_l}-\frac 1 n)(B^{t_l}-\tau_{\mathcal
M}(B^{t_l})+ \tau_{\mathcal M}(B^{t_l}))) \\
 &=  \tau_{\mathcal M}((e_{s_1s_1}-\frac 1 n)(B^{t_1}-\tau_{\mathcal
M}(B^{t_1}) )\ldots
  (e_{s_ls_l}-\frac 1 n)(B^{t_l}-\tau_{\mathcal
M}(B^{t_l}) ))\\
&\quad + \tau_{\mathcal M}(B^{t_1}) \tau_{\mathcal
M}((e_{s_1s_1}-\frac 1 n) (e_{s_2s_2}-\frac 1
n)(B^{t_2}-\tau_{\mathcal M}(B^{t_2}) )   \ldots
  (e_{s_ls_l}-\frac 1 n)(B^{t_l}-\tau_{\mathcal
M}(B^{t_l}) ))\\
&\quad + \ldots \\
 &\quad + \tau_{\mathcal M}(B^{t_1}) \tau_{\mathcal M}(B^{t_2})\ldots
 \tau_{\mathcal M}(B^{t_l})
 \tau_{\mathcal M}((e_{s_1s_1}-\frac 1 n)(e_{s_2s_2}-\frac 1
n)(e_{s_3s_3}-\frac 1 n)\ldots (e_{s_ls_l}-\frac 1 n))\\
&= 0 \\
&\quad + \tau_{\mathcal M}(X^{t_1}) \tau_{\mathcal
M}((e_{s_1s_1}-\frac 1 n) (e_{s_2s_2}-\frac 1
n)(X^{t_2}-\tau_{\mathcal M}(X^{t_2}) )  \ldots
  (e_{s_ls_l}-\frac 1 n)(X^{t_l}-\tau_{\mathcal M}(X^{t_l}) ))\\
&\quad + \ldots \\
&\quad + \tau_{\mathcal M}(X^{t_1}) \tau_{\mathcal M}(X^{t_2})\ldots
\tau_{\mathcal M}(X^{t_l})
 \tau_{\mathcal M}((e_{s_1s_1}-\frac 1 n)(e_{s_2s_2}-\frac 1
n)(e_{s_3s_3}-\frac 1 n)\ldots (e_{s_ls_l}-\frac 1 n))\\
& \qquad   \text{(because $\tau_{\mathcal
M}(B^i)=\tau_{\mathcal M} (X^i)$ for $1\le i\le k+1$ and $B$ (or $X$) is free from $\{e_{i_0i_0},e_{j_0j_0}\}$)}\\
 &  = \tau_{\mathcal M}((e_{s_1s_1}-\frac 1 n)X^{t_1}\ldots
  (e_{s_ls_l}-\frac 1 n)X^{t_l}),
  \end{aligned}$$
which proves $(iii)$  for $m=k+1$.

\vspace{0.3cm}

 \noindent As for $(iv)$, we know that
 $$\begin{aligned}
    \tau_{\mathcal M}(e_{i_0i_0}(Be_{j_0j_0})^{k+1}B) &=
    \frac 1 n \tau_{\mathcal N} (b_{i_0j_0}  (b_{j_0j_0})^{k}
    b_{j_0i_0})\\
    &= \frac 1 n \tau_{\mathcal N}   ((b_{j_0j_0})^{k}
     ) \tau_{\mathcal N} (
    b_{j_0i_0}b_{i_0j_0})
   \\
&= \frac 1 n \tau_{\mathcal N}  ((x_{j_0j_0})^{k}
     ) \tau_{\mathcal N}(
    x_{j_0i_0}x_{i_0j_0}) \qquad \text{(because of $(ii)$.)}\\
    &=
    \frac 1 n \tau_{\mathcal N} (x_{i_0j_0}  (x_{j_0j_0})^{k}
    x_{j_0i_0})\\
   &=\tau_{\mathcal M} (e_{i_0i_0}(Xe_{j_0j_0})^{k+1}X)
 \end{aligned}$$
 Now we are ready to prove $(v)$. We have
 $$
 \begin{aligned}
   \tau_{\mathcal M}& (e_{i_0i_0}(Be_{j_0j_0})^{k+1}B) \\&= \tau_{\mathcal M}((e_{i_0i_0}-\frac 1 n +\frac 1 n)B(e_{j_0j_0}-\frac 1 n
  +\frac 1 n)\ldots B(e_{j_0j_0}-\frac 1 n +\frac 1 n)B)\\
  &   =\tau_{\mathcal M}((e_{i_0i_0}-\frac 1 n  )B(e_{j_0j_0}-\frac 1 n  )\ldots B(e_{j_0j_0}-\frac 1 n
  )B) \qquad\quad \qquad \qquad \qquad \ (=0)\\
  &\quad + \frac 1 n \tau_{\mathcal M}( B^2(e_{j_0j_0}-\frac 1 n  )\ldots B(e_{j_0j_0}-\frac 1 n
  ))+ \ldots\\
 & \quad +\frac 1 {n^2} \tau_{\mathcal M}( B^3 (e_{j_0j_0}-\frac 1 n  ) \ldots B(e_{j_0j_0}-\frac 1 n
  )) +\ldots\\
  &\quad + \ldots\\
  &\quad + \frac 1 {n^{k+1}} \tau_{\mathcal M} ((e_{i_0i_0}-\frac 1 n  )
  B^{k+2}) + \frac 1 {n^{k+1}} \tau_{\mathcal M} ((e_{j_0j_0}-\frac 1 n  )
  B^{k+2})+ \ldots \qquad\quad  (=0)\\
  &\quad + \frac 1 {n^{k+2}}\tau_{\mathcal M} (B^{k+2})
 \end{aligned}
 $$
On the other hand, similar computation shows
$$
 \begin{aligned}
   \tau_{\mathcal M}& (e_{i_0i_0}(Xe_{j_0j_0})^{k+1}X) \\ &=   \tau_{\mathcal N}((e_{i_0i_0}-\frac 1 n +\frac 1 n)X(e_{j_0j_0}-\frac 1 n
  +\frac 1 n)\ldots X(e_{j_0j_0}-\frac 1 n +\frac 1 n)X)\\
  &   =\tau_{\mathcal M}((e_{i_0i_0}-\frac 1 n  )X(e_{j_0j_0}-\frac 1 n  )\ldots X(e_{j_0j_0}-\frac 1 n
  )X) \qquad\quad \qquad \qquad \qquad \ (=0)\\
  &\quad + \frac 1 n \tau_{\mathcal M}( X^2(e_{j_0j_0}-\frac 1 n  )\ldots X(e_{j_0j_0}-\frac 1 n
  ))+ \ldots\\
 & \quad +\frac 1 {n^2} \tau_{\mathcal M}( X^3 (e_{j_0j_0}-\frac 1 n  ) \ldots X(e_{j_0j_0}-\frac 1 n
  )) +\ldots\\
  &\quad + \ldots\\
  &\quad + \frac 1 {n^{k+1}} \tau_{\mathcal M} ((e_{i_0i_0}-\frac 1 n  )
  X^{k+2}) + \frac 1 {n^{k+1}} \tau_{\mathcal M} ((e_{j_0j_0}-\frac 1 n  )
  X^{k+2})+ \ldots \qquad\quad  (=0)\\
  &\quad + \frac 1 {n^{k+2}}\tau_{\mathcal M} (X^{k+2})
  \\
 \end{aligned}
 $$
Combining with $(iii)$ and $(iv)$, we know that $\tau_{\mathcal
M}(B^{k+2}) =\tau_{\mathcal M}(X^{k+2}) $. i.e. $(v)$ is proved for
$m=k+1$.
 Q.E.D

\vspace{0.5cm}

  Note that   Voiculescu's semicircular matrix  always
satisfies all conditions in Proposition 3.1. Now  we have
\begin{Th}
Suppose that $B=\sum_{1\le i,j\le n} b_{ij}\otimes e_{ij}$,is a
self-adjoint element in $\mathcal M$ (with $b_{ij}$'s   in $\mathcal
N$) such that, for some $i_0, j_0, \ 1\le i_0 \ne j_0\le n$,
\begin{enumerate}
\item $B$ and $\{I_{\N} \otimes e_{i_0i_0},I_{\N} \otimes e_{j_0j_0}\}$ are free  with respect to
 $\tau_{\mathcal M}$;
\item  $b_{ i_0j_0}^*b_{i_0j_0}$ and  $b_{j_0j_0}$ are free in $\mathcal N$ are free  with respect to
 $\tau_{\mathcal N}$.
 \end{enumerate}
 Then $  B $ is a semicircular element of $(\tau_{\mathcal M}(B),
 \|B\|_2)$.
\end{Th}

\noindent {\bf Proof:} Let $X$ be a Voiculescu's semicircular
matrix, and $\tilde B= \frac 1 {\|B-\tau_{\mathcal M}(B)\|_2}
(B-\tau_{\mathcal M}(B))$. Using Proposition 3.1, we have that
$\tau_{\M}(\tilde B^m) = \tau_{\M}(X^m)$ for every $m\in \Bbb N$.
Hence $B$ is a semicircular element of $(\tau_{\mathcal M}(B),
 \|B\|_2)$. Q.E.D.

\vspace{0.3cm} \noindent By switching the order of $i$ and $j$, same
arguments as Corollary 3.1 will show the following.
\begin{Th}
Suppose that $B=\sum_{1\le i,j\le n} b_{ij}\otimes e_{ij}$,is a
self-adjoint elements in $\mathcal M$ (with $b_{ij}$'s   in
$\mathcal N$) such that, for some $i_0, j_0, \ 1\le i_0 \ne j_0\le
n$,
\begin{enumerate}
\item $B$ and $\{I_{\N} \otimes e_{i_0i_0},I_{\N} \otimes  e_{j_0j_0}\}$  are free  with respect to
 $\tau_{\mathcal M}$;
\item  $b_{i_0j_0} b_{i_0j_0}^*$ and  $b_{i_0i_0}$ are free  with respect to
 $\tau_{\mathcal N}$.
 \end{enumerate}
 Then $  B $ is a semicircular element of $(\tau_{\mathcal M}(B),
 \|B\|_2)$.
\end{Th}

\vspace{0.5cm}

  Our next result describes the uniqueness of Voiculescu's
   semicircular matrix. Let $\mathcal D_n$ be the abelian von Neumann
subalgebra generated by $\{I_{\N} \otimes e_{ii}\}_{1\le i\le n}$ in
$\mathcal M$.
\begin{Prop}
Suppose that $B=\sum_{1\le i,j\le n} b_{ij}\otimes e_{ij}$  and
$X=\sum_{1\le i,j\le n} x_{ij}\otimes e_{ij}$ are two self-adjoint
elements  in $\mathcal M$ (with $b_{ij}$'s and $x_{ij}$'s  in
$\mathcal N$), such that
\begin{enumerate}
  \item   $B$ and   $\mathcal D_n$ are free  with respect to
 $\tau_{\mathcal M}$;
  \item $X$  and   $\mathcal D_n$ are free  with respect to
 $\tau_{\mathcal M}$;
  \item $\tau_{\mathcal M}(B^m) = \tau_{\mathcal M}(X^m)$  for all
  $m\ge 1$;
  \item $b_{in}, x_{in}$ are positive elements in $\mathcal N$ for
  all $1\le i\le n-1$;
  \item   The spectrum of $x_{in}$ has no atom as the operator in $\mathcal N$ for
  all $1\le i\le n-1$;
\end{enumerate}
Then, $\{b_{ij}\}_{1 \le i,j\le n}$  and $\{x_{ij}\}_{1 \le i,j\le
n}$ have the  identical $*$-joint distribution in $\N$  with respect
to the tracial state $\tau_{\mathcal N}$.

\end{Prop}

\noindent {\bf Proof: } Note that    $B$ (or $X$) is free from
$\mathcal D_n$ in $\mathcal M$ and $\tau_{\mathcal M}(B^m) =
\tau_{\mathcal M}(X^m)$  for all
  $m\ge 1$. We have
$$\begin{aligned}
\tau_{\mathcal M}((I_{\N}\otimes e_{i_1i_1})&B(I_{\N}\otimes
e_{i_2i_2})B\cdots (I_{\N}\otimes e_{i_mi_m})B)\\ &=\tau_{\mathcal
M} ((I_{\N}\otimes e_{i_1i_1})X(I_{\N}\otimes e_{i_2i_2})X\cdots
(I_{\N}\otimes e_{i_mi_m})X),\end{aligned}
$$ or
\begin{equation}
\tau_{\mathcal N}(b_{i_1i_2}b_{i_2i_3} \cdots b_{i_mi_1}
)=\tau_{\mathcal N} (x_{i_1i_2}x_{i_2i_3} \cdots x_{i_mi_1}),
\end{equation} for all $i_1,\ldots, i_m \in \{1,\ldots, n\}$. \ In
particular, we have
\begin{align}\tau_{\mathcal N}((b_{in}b_{ni})^m) =\tau_{\mathcal
N}((x_{in}x_{ni})^m),\quad \text {for all $m\ge 0$ }. \end{align}

\noindent Since $x_{in}$ is a positive element whose spectrum has no
atom as the operator in $\mathcal N$ for
  all $1\le i\le n-1$,   by functional
calculus, we know that
 \begin{align} I_{\mathcal N}\in \{x_{in}\}''&=\text{the von Neumann subalgebra generated by $x_{in}$  in $\mathcal N$.}\notag\\
&=\overline{span}^{SOT}\{x_{in}^{2t}: 1\le t \in \Bbb N\} \tag{$*$} \\
& =\overline{span}^{SOT}\{x_{in}^{2t-1}: 1\le t \in \Bbb
N\}\notag\end{align}
 By (3), we have
\begin{align} I_{\mathcal N}\in  \{b_{in}\}''&=\text{the von Neumann subalgebra generated by $b_{in}$  in $\mathcal N$.}\notag\\
&=\overline{span}^{SOT}\{b_{in}^{2t}: 1\le t \in \Bbb N\} \tag{$**$} \\
& =\overline{span}^{SOT}\{b_{in}^{2t-1}: 1\le t \in \Bbb
N\}\notag\end{align}

 \noindent To prove the proposition, we will show  the following lemma first.
\begin{lemma}Following the notions as above, we have that, for any
$1\le i_1, \ldots, i_{2m} \le n $,
\begin{align}
\tau_{\mathcal N}(b_{i_1i_2}b_{i_3i_4}\cdots b_{i_{2m-1}i_{2m}})
   =\tau_{\mathcal N}(x_{i_1i_2}x_{i_3i_4}\cdots x_{i_{2m-1}i_{2m}})
\end{align}
\end{lemma}
\noindent{\bf Proof of Lemma: } For    $1\le i_1, \ldots, i_{2m} \le
n$, let $i_{2m+1}$ be $ i_1$ and  $$\mathcal S[i_1, i_2,\ldots,
i_{2m}] = \{i_{2j}\ | \ i_{2j} \ne i_{2j+1}, \ 1\le j\le m\} \subset
\{i_2, i_4, \ldots, i_{2m} \}.$$  Denote by $l$ the cardinality of
the set $\mathcal S[i_1, i_2,\ldots, i_{2m}]$.   There are two cases
to be considered,  $(i) \ m=1$ and $(ii) \ m\ge 2$.

\vspace{0.3cm}

\noindent $(i)$: \noindent When $m=1$ and $l=0$, (4) follows
directly from (2).  When  $m=l=1$,  we need to show that
$\tau_{\mathcal N}(b_{i_1i_2})= \tau_{\mathcal N}(x_{i_1i_2})$ for
$i_1\ne i_2$. If one of $i_1$ or $i_2$ is equal to $n$, then
$\tau_{\mathcal N}(b_{i_1i_2}b_{i_2i_1})= \tau_{\mathcal
N}(b_{i_1i_2}b_{i_1i_2}^*)=\tau_{\mathcal
N}(x_{i_1i_2}x_{i_2i_1})=\tau_{\mathcal N}(x_{i_1i_2}x_{i_1i_2}^*)$
by (3). Therefore $\tau_{\mathcal N}(b_{i_1i_2})= \tau_{\mathcal
N}(x_{i_1i_2})$ because both $b_{i_1i_2}$ and $x_{i_1i_2}$ are
positive. If none of $i_1$ and $i_2$ is equal to $n$, it induces
from (2) that, for all $t_1, t_2$ in $\Bbb N \cup \{0\}$,
$$ \tau_{\mathcal N}(b_{ni_1}(b_{i_1n}b_{ni_1})^{t_1}b_{i_1i_2}(b_{i_2n}b_{ni_2})^{t_2}b_{i_2n})
        =\tau_{\mathcal N}
        (x_{ni_1}(x_{i_1n}x_{ni_1})^{t_1}x_{i_1i_2}(x_{i_2n}x_{ni_2})^{t_2}x_{i_2n}).
 $$
 or (note that both $b_{i_1n}$ and $x_{i_1n}$ are
positive.)
 $$ \tau_{\mathcal N}( (b_{i_1n} )^{2t_1+1}b_{i_1i_2}(b_{i_2n} )^{2t_2+1} )
        =\tau_{\mathcal N}
        ((x_{i_1n} )^{2t_1+1}x_{i_1i_2}(x_{i_2n} )^{2t_2+1}), \ \ \forall t_1,t_2 \in \Bbb N\cup \{0\}.
 $$
From $(*), (**)$ it follows that
  $$  \tau_{\mathcal
N}(b_{i_1i_2})= \tau_{\mathcal N}(x_{i_1i_2})  $$

 \vspace{0.3cm}

\noindent $(ii)$: When $m\ge 2$, we will use induction on $l$ now.

 Since $(4)$ holds when $l=0$, we might assume that   $(4)$ holds for $l\le k$.
Consider the case when $l=k+1$. We need to show that
\begin{align}
\tau_{\mathcal N}(b_{i_1i_2}b_{i_3i_4}\cdots b_{i_{2m-1}i_{2m}})
   =\tau_{\mathcal N}(x_{i_1i_2}x_{i_3i_4}\cdots x_{i_{2m-1}i_{2m}}),
\end{align}
 when the cardinality of the set $\mathcal S[i_1, i_2,\ldots,
i_{2m}]$
 is equal to $k+1$. Since $\tau_{\mathcal N}$ is a
tracial state, we can assume that $i_2\ne i_3.$ There are three
cases we have to consider: (a) $i_2=n,\ i_3\ne n$, (b) $i_2\ne n, \
i_3=n$, (c) $i_2\ne n, \ i_3\ne n.$

 \vspace{0.3cm}

  (a): Assume that $i_2=n,\ i_3\ne n$. Since cardinality of $\mathcal S[i_1, i_2,i_3,i_4\ldots, i_{2m}]$ is equal to $k+1$, the
cardinality of
$$
   \mathcal S[i_1, i_2,\overbrace {i_2, i_3,i_3,i_2},\ldots, \overbrace {i_2, i_3,i_3,i_2}, \overbrace { i_2,i_3}, i_3, i_4,\ldots, i_{2m}]
$$ is equal to $k$ (here $i_2=n$).
  By hypothesis, we have that
$$
\begin{aligned}
\tau_{\mathcal N}(b_{i_1n}(b_{ni_3})^{2t-1}b_{i_3i_4}\cdots
b_{i_{2m-1}i_{2m}})
   =\tau_{\mathcal N}(x_{i_1n}(x_{ni_3})^{2t-1}x_{i_3i_4}\cdots
   x_{i_{2m-1}i_{2m}}), \forall t\in \Bbb N.\end{aligned}
$$
From $(*), (**)$, we have that  (5) holds for this family of
$\{i_1,\ldots, i_{2m}\}.$

 \vspace{0.3cm}

(b): The proof of the case when $i_2\ne n, \ i_3=n$ is similar as
(a).

 \vspace{0.3cm}

(c): Assume that  $i_2\ne n, \ i_3\ne n.$ Since cardinality of
$\mathcal S[i_1, i_2,i_3,i_4\ldots, i_{2m}]$ is equal to $k+1$, the
cardinality of
$$
   \mathcal S[i_1, i_2,\overbrace{i_2, n,n,i_2},\ldots,\overbrace{i_2, n,n,i_2},\overbrace{i_2, n},\overbrace{ n,i_3},\overbrace{i_3, n,n,i_3}, \ldots,\overbrace{i_3, n,n,i_3},
    i_3, i_4,\ldots, i_{2m}]
$$ is equal to $k$ (here both $i_2, i_3$ are not equal to $n$).
By hypothesis, we have that
$$
\begin{aligned}
 \tau_{\mathcal
N}(b_{i_1i_2}(b_{i_2n})^{2t_1-1}&(b_{ni_3})^{2t_2-1}b_{i_3i_4}\cdots
b_{i_{2m-1}i_{2m}})\\&
   =\tau_{\mathcal N}(x_{i_1i_2}(x_{i_2n})^{2t_1-1}(x_{i_3n})^{2t_2-1}x_{i_3i_4}\cdots
   x_{i_{2m-1}i_{2m}}), \ \ \  \forall t_1, t_2\in \Bbb N.\end{aligned}
$$
From $(*), (**)$, we have that  (5) holds for this family of
$\{i_1,\ldots, i_{2m}\}.$

 \vspace{0.3cm}

Hence the induction is completed, and the lemma is proved.  Q.E.D.

\vspace{0.5cm}

\noindent {\bf Continue the proof of the Proposition 3.2:} Note that
$B$ and $X$ are two self-adjoint elements in $\mathcal M$. Thus
$b_{ij}=b_{ji}^*$ and $x_{ij}=x_{ji}^*$. From the preceding Lemma,
we know that the *-joint distribution of the family of elements $\{
a_{ij}, 1\le i,j\le n \}$ is identical to the *-joint  distribution
of the family of elements $\{x_{ij},  1\le i,j\le n\}$ in $\N$ with
respect to $\tau_{\N}$. Q.E.D.

\vspace{0.5cm}

Now we can prove our main theorems in this section.
\begin{Th}
  Suppose $B= \sum_{1\le i,j\le n} b_{ij}\otimes e_{ij}$ is a
  semicircular element of $(0,1)$ in $\mathcal M$. If $B$ and $\mathcal
  D_n$  are free with respect to $\tau_{\mathcal M}$,  then there is a family of unitary elements
  $\{u_1, \ldots, u_n\}$ in $\mathcal N$ such that $U^*BU$ is a
  Voiculescu's semicircular matrix, where $U=\sum_{1\le
  i\le n} u_i\otimes e_{ii}$.
\end{Th}

\noindent {\bf Proof: } Let $ u_ih_i$ be the polar decomposition of
$b_{in}$ in $\N$ for $1\le i \le n-1$. Let $u_n= I_{\N}$,
$U=\sum_{1\le
  i\le n} u_i\otimes e_{ii}$ and $\tilde B =U^*BU =\sum_{1\le i,j\le n} u_i^*b_{ij}u_j\otimes e_{ij} $. Let $X$ be
a
  Voiculescu's semicircular matrix. It is easy to check that $\tilde
  B$ and $X$ satisfy the conditions in Proposition 3.2. Hence $\{
  u_i^*b_{ij}u_j\}_{1\le i,j \le n}$ and $\{
  x_{ij}\}_{1\le i,j \le n}$ have identical *-joint distributions in
  $\N$.  It follows that $\tilde B$ is also a Voiculescu's
  semicircular matrix. Q.E.D.

\vspace{0.3cm}

The proofs of Proposition 3.2 and Theorem 3.3 can be easily extended
to the case of $m$-tuple of semicircular matrices, $B_1, \ldots,
B_m$ when $\{B_1, \ldots, B_m, \mathcal D_n\} $ are free with
respect to $\tau_{\M}$. We present the following theorem whose proof
is skipped.
\begin{Th}
  Suppose $B_{k}= \sum_{1\le i,j\le n} b_{ij}^{(k)}\otimes e_{ij}$  for $1\le k \le m$ is
  a sequence of
   semicircular elements of $(0,1)$ in $\mathcal M$. If $\{B_{1},\ldots, B_{m}, \mathcal D_n\}$ are free
    in $\mathcal M$, then there is a family of unitary elements
  $\{u_1, \ldots, u_n\}$ in $\mathcal N$ such that $\{U^*B_{k}U\}_{1\le k \le m}$ is a
  standard family of
  Voiculescu's  semicircular matrices, where $U=\sum_{1\le
  i\le n} u_i\otimes e_{ii}\in \M$.
\end{Th}
\noindent{\bf Remark: } Theorem 3.3 can be viewed as the inverse
statement of (Voiculescu's) Theorem 2.1.

\begin{Col}
Suppose $B_{k}= \sum_{1\le i,j\le n} b_{ij}^{(k)}\otimes e_{ij}$ for
$1\le k \le m$ is
  a sequence of semicircular elements of $(0,1)$
   in $\mathcal M$. Then  $\{B_{1},\ldots, B_{m}, \mathcal D_n\}$  are free in $\mathcal M$ if and only if there is a family of unitary elements
  $\{u_1, \ldots, u_n\}$ in $\mathcal N$ such that $\{U^*B_{k}U\}_{1\le k \le m}$ is
  a
  standard family of
  Voiculescu's  semicircular matrices, where $U=\sum_{1\le
  i\le n} u_i\otimes e_{ii}\in \M$.
\end{Col}

 \vspace{0.3cm}
 \noindent Combining with Theorem 3.1, we have the following
characterization of semicircular elements in $\mathcal M$.
 \begin{Th}
    Suppose $B= \sum_{1\le i,j\le n} b_{ij}\otimes e_{ij}$ is a
  self-adjoint element in $\mathcal M$  and free from $\mathcal
  D_n$ with respect to $\tau_{\M}$. Then the following are equivalent.
  \begin{enumerate}[(i)]
      \item $B$ is a semicircular element in $\mathcal M$.
      \item There is a family of unitary elements
  $\{u_1, \ldots, u_n\}$ in $\mathcal N$ such that $U^*BU$ is a
  Voiculescu's   semicircular matrix, where $U=\sum_{1\le
  i\le n} u_i\otimes e_{ii}$.
  \item For all $i , j , \ 1\le i  \ne j \le
n$, $b_{ji}^*b_{ij} $ and  $b_{jj}$ are free in $\mathcal N$ with
respect to the tracial state $\tau_{\mathcal N}$.
\item There are some  $i_0, j_0, \ 1\le i_0 \ne j_0\le
n$, $b_{i_0j_0}^*b_{i_0j_0} $ and  $b_{j_0j_0}$ are free in
$\mathcal N$ with respect to the tracial state $\tau_{\mathcal N}$.
\item For all $i , j , \ 1\le i  \ne j \le
n$, $b_{ij} b_{ij}^*$ and  $b_{ii}$ are free in $\mathcal N$ with
respect to the tracial state $\tau_{\mathcal N}$.
\item There are some  $i_0, j_0, \ 1\le i_0 \ne j_0\le
n$, $b_{ i_0j_0} b_{i_0j_0}^*$ and  $b_{i_0i_0}$ are free in
$\mathcal N$ with respect to the tracial state $\tau_{\mathcal N}$.
  \end{enumerate}
 \end{Th}
\noindent {\bf Proof: } $(i) \Rightarrow (ii)$ is from Theorem 3.3.
$(ii) \Rightarrow (iii)$ is from Definition 2.5. $(iii)\Rightarrow
(iv)$ is trivial. $(iv) \Rightarrow (i)$ is from Theorem 3.1. $(ii)
\Rightarrow (v)$ is from Definition 2.5. $(v)\Rightarrow (vi)$ is
trivial. $(vi) \Rightarrow (i)$ is from Theorem 3.2. Q.E.D.

 \vspace{0.3cm}
\noindent As a corollary of Theorem 3.4, we prove the following
result of ``free compression" in [7].
\begin{Col} Suppose $A_{k}= \sum_{1\le i,j\le n} a_{ij }^{(k)}\otimes e_{ij}$ for $1\le k\le m$ is a
  sequence of self-adjoint elements in $\mathcal M$. If $\{A_{1},\ldots, A_{m}, \mathcal D_n\}_{1\le k\le m}$ are free in $\mathcal M$, then
\begin{enumerate}[(i)]
\item
  $\{ a_{ii }^{(1)},\ldots,  a_{ii }^{(m)}\}$ is a free
  family in $\mathcal N$ for each $1\le i\le n$.
  \item   $ \{a_{ij_1 }^{(1)}a_{ij_1 }^{(1)*},\ \ldots, \  a_{ij_m }^{(m)} a_{ij_m }^{(m)*}\}$ is a free
  family in $\mathcal N$ for each $1\le i, j_1, \ldots, j_m\le n$.\end{enumerate}
\end{Col}

\noindent {\bf Proof: } We can assume that there exists a family of
semicircular elements $\{X_{1}, \ldots, X_{m}\}$ in $\M$ such that
$\{X_{1},\ldots, X_{m}, \mathcal D_n\}_{1\le k\le m}$ are free  in
$\mathcal M$ and each $A_{k}$ is contained in the von Neumann
subalgebra generated by $X_{k}$ in $\M$ for $1 \le k \le m$. It
follows from Theorem 3.4 that there exists a family of unitary
elements
  $\{u_1, \ldots, u_n\}$ in $\mathcal N$ such that $\{U^*X_{k}U\}_{1\le k \le m}$ is a
  standard family of
  Voiculescu's  semicircular matrices, where $U=\sum_{1\le
  i\le n} u_i\otimes e_{ii}$. Assume that $X_{k} = \sum_{1\le i,j\le n} x_{ij }^{(k)}\otimes
  e_{ij}$. Therefore, from the definition of   standard family of
  Voiculescu's semicircular matrices, we have
$\{Re\  u^*_ix_{ i,j }^{(k)}u_j, Im\  u^*_ix_{ i,j }^{(k)}u_j \  |\
1 \le i\le j \le n, 1\le k\le m\} $ is a family of free elements in
$\mathcal N$ with respect to $\tau_{\mathcal N}$.

 Note that each $A_{k}$ can be
approximated by the polynomials of $X_{ k }$ in 2-norm $\|\ \cdot
\|_2$. So $u_i^*a_{ii}^{(k)}u_i $ can be approximated by the
polynomials of $\{Re\ u^*_sx_{ s,t }^{(k)}u_t, Im\  u^*_sx_{ s,t
}^{(k)}u_t \ |\ 1 \le s, t \le n\} $ in 2-norm $\|\ \cdot \|_2$.
Since  $\{Re\ u^*_ix_{ i,j }^{(k)}u_j, Im\  u^*_ix_{ i,j }^{(k)}u_j
\ |\ 1 \le i\le j \le n, 1\le k\le m\} $ is a family of free
elements in $\mathcal N$ with respect to $\tau_{\mathcal N}$, we
have that
  $\{ u^*_ia_{ii }^{(1)}u_i ,\ldots,  u^*_ia_{ii }^{(m)}u_i $ \} is a free
  family in $\mathcal N$ for each $1\le i\le n$. Hence $\{  a_{ii }^{(1)}  ,\ldots,  a_{ii }^{(m)}  $ \} is a free
  family in $\mathcal N$ for each $1\le i\le n$.  Similarly, we also have
  $ \{a_{ij_1 }^{(1)}a_{ij_1 }^{(1)*},\ \ldots, \  a_{ij_m }^{(m)} a_{ij_m }^{(m)*}\}$ is a free
  family in $\mathcal N$ for each $1\le i, j_1, \ldots, j_m\le n$.

\vspace{2cm}

\section{Matricial distance between two free
semicircular elements.}

Let $\mathcal M$ be a factor of type II$_1$ with the traical state
$\tau_{\mathcal M}$.  
Assume that  $\mathcal M \cong \mathcal N\otimes \mathcal M_n$ for
some positive integer $n$ and a type II$_1$ subfactor $\mathcal N$
with the tracial state  $\tau_{\N}$. Let $\{e_{ij}\}_{i,j=1}^n$ be
the canonical system of matrix units of $\mathcal M_n$ and $\mathcal
D_n$
  be the subalgebra generated by   $\{I_{\N}\otimes e_{ii}\}_{1\le i \le
 n}$ in $\mathcal M$.
 Let $E_n$ be the conditional expectation
  from $\mathcal M$ onto $I_{\N}\otimes M_n$.

\begin{Prop}
  If $B$ is a
  semicircular element of (0,1) in $\mathcal M$, which is free from $\mathcal
  D_n$    with respect to the trace $\tau_{\mathcal
  M}$, then
  $$ \| E_n(B)\|_2 \le \frac {7} {n^{1/8}}.  $$
\end{Prop}

\noindent The proof of the proposition will be given in later
subsection. We need to prove a few lemmas first.

\subsection{Definition of ${\tilde{\mb x}}$:  perturbation of a semicircular element ${\mb x}$ }
Let $\mathcal N_1$ be a type II$_1$ von Neumann algebra with a
traical state $\tau$. Let $u$ be a Haar unitary element in $\N_1$
and  $\mathcal A_1$ be the diffused abelian von Neumann subalgebra
generated by $u$ in $\mathcal N_1$. There is a $*$-isomorphism
$\psi$ from $\A_1$ onto $L^{\infty}([0,1],m)$ such that $\psi(u)=
e^{2\pi i\cdot t}$, where $m$ is the Lebesgue measure on $[0,1]$.
Let $$ s(t) =\int_{-r}^t \frac 2 {\pi r^2}\sqrt {r^2-t_1^2}dt_1, \ \
\qquad \text{for } t \in [-r,r].
 $$ So,
 $$ \begin{aligned}
   \frac {ds}{dt}  =\frac 2{\pi r^2} \sqrt {r^2-t^2}; \qquad \quad \  \frac {dt}{ds}  =\frac  {\pi r^2}  2 \frac 1 {\sqrt
  {r^2-t^2}},
\end{aligned} $$and  $$ \begin{aligned}
   s(-r) =0; \qquad \quad  \qquad \quad   \quad  s(r)=1.
\end{aligned}
  $$
Let $t= g(s)$  be the inverse function of $s=s(t)$. Then $g(s) \in
L^\infty([0,1],m)$, $g(0)=-r,$ $g(1)=r$, and $\frac {dt}{ds} =
g'(s)$. Since
$$ \begin{aligned}
  &  \int_{0}^1 g(s)^m ds = \int_{-r}^r g(s(t))^m \frac
2{\pi r^2} \sqrt {r^2-t^2} dt = \frac 2{\pi r^2} \int_{-r}^r  t^m
\sqrt {r^2-t^2} dt,
\end{aligned}
$$ we know that  $\psi^{-1}(g(s))$ defines an operator $x$,  a semicircular element of $(0,r)$,  in $\A_1$.
  Let
$$
 f(s) =\left \{ \begin{aligned}
   & \frac {g(r)}{r^2}s^2 \qquad \quad \qquad  \qquad   0\le s \le r\\
   &  g(s) \qquad \quad \qquad\qquad \quad  r< s \le 1-r\\
   &   \frac {g(1-r)}{r^2} (1-s)^2  \qquad \quad 1-r \le s \le 1.
 \end{aligned}\right .
$$

\begin{Def}
\noindent Define $\tilde x$, the perturbation of a semicircular
element $x$, to be $\psi^{-1}(f(s))$, the 
corresponding element of $f(s)$ in $\mathcal A_1$.
\end{Def}
\noindent
\begin{lemma} Assume that  $x$ is  a semicircular element of $(0,r)$ in the
 abelian von Neumann subalgebra $\A_1$ generated by the Haar unitary $u$. And
$\tilde x$ is defined as   in Def. 4.1. Then we have
$$
\sum_{k\ne 0}\left |\tau( \tilde  x u^{ k}   )\right | \le
    \sqrt 5 r^{\frac 1 2}
$$
\end{lemma}
\noindent {\bf Proof:} From the construction of $\tilde
x$, we have
$$
   \begin{aligned}
      \tau(\tilde x u^k)  = \int_0^1 f(s) e^{2\pi i ks}ds
       = - \frac 1 {2\pi i k} \int_0^1 f'(s)  e^{2\pi i ks} ds
   \end{aligned}
$$
Hence
$$\begin{aligned}
  \left (\sum_{k\ne 0}\left |\tau( \tilde  xu^{ k}   )
  \right | \right )^2 & \le  \left (\sum_{k\ne 0 }\left  |\frac 1 {2\pi i k}\right | \cdot \left | \int_0^1 f'(s)  e^{2\pi i
  ks}
  ds \right | \right )^2\\
  &\le \left (\sum_{k\ne 0} \left |\frac 1 {2\pi i k} \right |^2 \right )\left  (\sum_{k\ne 0} \left | \int_0^1 f'(s)  e^{2\pi i
  ks}
  ds \right |^2 \right ) \le   \int_0^1 |f'(s)|^2
ds.   \end{aligned}$$ Because
$$
  \begin{aligned}
    \int_0^r  |f'(s)|^2 ds & =\int_0^r  \frac {|g(r)|^2}{r^4}(2s)^2
    ds = \frac 4 3  \frac {|g(r)|^2}{r} \le \frac 4 3  r  \qquad (\text{ since $|g(r)|\le r$})\\
   \int_{1-r}^1  |f'(s)|^2 ds &\le   \frac 4 3 r \qquad\qquad \qquad (\text{ similar as the preceding one})
   \\
 \int_r^{1-r}
|f'(s)|^2 ds  & =  \int_r^{1-r} |g'(s)|^2 ds =
\int_{g(r)}^{g(1-r)}\left | \frac {dt}{ds}  \right |^2 \frac {ds}
{dt} dt\\
&= \int_{g(r)}^{g(1-r)} \left ( \frac {\pi r^2}{2} \right )^2 \frac
{1}{r^2-t^2} \cdot \frac 2 {\pi r^2} \sqrt{r^2-t^2} dt\\
& = \frac   {\pi r^2}2 \int_{g(r)}^{g(1-r)} \frac {d(\frac t
r)}{\sqrt {1-(\frac t r)^2}} \\
&\le \frac  {\pi r^2} 2\int_{-1}^{1} \frac {d t  }{\sqrt {1-  t ^2
}} \le r, \qquad\qquad (\text{ when $r$ is small enough.})
  \end{aligned}
$$
we have
$$ \begin{aligned}
\int_0^1 |f'(s)|^2 ds & = \int_0^r  |f'(s)|^2 ds+  \int_r^{1-r}
|f'(s)|^2 ds + \int_{1-r}^1  |f'(s)|^2 ds, \\
&\le \frac 4 3 r + \frac 4 3 r+ r  \le 5r;\qquad\\
    (\sum_{k\ne 0}\left |\tau( \tilde  xu^{ k}   )\right |)^2
 &\le\int_0^1 |f'(s)|^2 ds
  \le  5 r;\qquad\\
       \sum_{k\ne 0}\left |\tau( \tilde  xu^{ k}   )\right | &\le
    \sqrt 5 r^{\frac 1 2}\end{aligned}
$$ Q.E.D.

\vspace{0.5cm}

\noindent Using the preceding notations, we have the following
inequality.
 \begin{lemma} For any two elements
$\{w_1, w_2\}$ in $\mathcal N_1$ with $||w_i||\le 1,$ we have
$$
  \begin{aligned}
          |\tau(w_1xw_2 )|^2   \le  |\tau(w_1\tilde x
    w_2 )|^2 +6r^{\frac 5 2}.
  \end{aligned}
$$
\end{lemma}

\noindent {\bf Proof: } By the definition of $\tilde x$,
 we get,
$$
\begin{aligned}
   \|x-\tilde x\|_2^2  &= \int_0^1|g(s)-f(s)|^2ds\\
   &\le \int_0^r|g(s)-\frac {g(r)}{r^2}s^2|^2ds + \int_{1-r}^1|g(s)-\frac {g(1-r)}{r^2} (1-s)^2
   |^2ds\\
   &\le \int_0^r|g(s) |^2ds + \int_{1-r}^1|g(s)  |^2ds\\
  \text {(as } g(r)\le 0 ,\  g(s)\le 0 & \text{ when }  0\le s\le r \text{
  and }   \ g(1-r)\ge 0 , \
    g(s)\ge 0  \text{ when } 1-r\le s\le 1)\\
    &\le r^2\cdot r + r^2\cdot r =2r^3.
\end{aligned}
$$
Therefore, for any two elements $w_1, w_2$ in $\mathcal N_1$ with
$\|w_1||, ||w_2|| \le 1$, we have
$$
 \begin{aligned}
   |\tau(w_1xw_2 )-\tau(w_1\tilde x w_2 )|^2
   \le || x-\tilde x||^2_2\le 2 r^3.
 \end{aligned}
$$
And $$
       \begin{aligned}
         & \left ||\tau(w_1xw_2 )|^2-|\tau(w_1\tilde x w_2 )|^2 \right |\\
             &  \quad \le |\tau(w_1xw_2 ) - \tau(w_1\tilde x w_2 )| \cdot (|\tau(w_1xw_2 )|+| \tau(w_1  x
             w_2 ) -(\tau(w_1xw_2 )-\tau(w_1\tilde x w_2 )) |) \\
             &\quad \le \sqrt 2 r^{\frac 3 2}(||x||_2+||x||_2 +\sqrt 2 r^{\frac 3
             2}) < 6r^{\frac 52}
       \end{aligned}
    $$
It follows that,
$$
  \begin{aligned}
          |\tau(w_1xw_2 )|^2
    \le  |\tau(w_1\tilde x
    w_2 )|^2 +6r^{\frac 5 2}.
  \end{aligned}
$$


\subsection{Some analysis on free group factors}

Let $\Sigma$ be an index set, $F(\Sigma)$ be the free group with the
standard generators $\{g_{\alpha}\}_{\alpha\in \Sigma}$. Let
$\lambda$ be the left regular representation of $F(\Sigma)$ and
$L(F(\Sigma))$ be the free group factor associated with the group
$F(\Sigma)$ with standard generators
$\{\lambda(g_{\alpha})\}_{\alpha \in \Sigma}$.

 Fix some index  $\alpha $ in $\Sigma$. We are going to group the
elements of $F(\Sigma)$ into following sets.Let
$$
\begin{aligned}
 ES  & = \{ w \in F(\Sigma) \ | \ w  \ \text { is a reduced word ending with } g_{\alpha}^{m}  \text{ such that
 } m  \ne 0 \}\\
 SS  & = \{ w \in F(\Sigma) \ | \ w  \ \text { is a reduced word starting with } g_{\alpha}^n \text{  such that
 } n \ne 0 \}\\
 ET &= F(\Sigma)\setminus ES, \qquad \qquad \\
 ST  &= F(\Sigma)\setminus ET.
\end{aligned}
$$
 \noindent Note that every element $g$ in $ES$ can be
expressed as $ag_{\alpha}^{m} $ for some reduced word $a$ (not
ending with $g_{\alpha} $) and nonzero integers $m$; and every
element $h$ in $SS$ can be expressed as $g_{\alpha }^{n} b$ for some
reduced word $b$  (not starting with $g_{\alpha} $) and nonzero
integers $n$.
 Then for every
$w$ in $L(F(\Sigma))$, we   get the expression of $w$ as
$$
\begin{aligned}
  w &=  \sum_{g\in ET\cup ES  }w(g) u_g = \sum_{   ag_{\alpha}^m\in ES  }\mathfrak{E}(a, m) \lambda(a)\lambda(g_{\alpha})^m
  + \sum_{  a \in ET  }\mathfrak{E}(a,0 ) \lambda(a)
   \\
  w &=  \sum_{h\in ST \cup SS }w(h) u_h= \sum_{    g_{\alpha}^n b\in SS  }\mathfrak{S}(a, m)
  \lambda(g_{\alpha})^n\lambda(b)
  + \sum_{ b \in ST  }\mathfrak{S}(b,0 )
  \lambda(b)
\end{aligned}
$$ where $w(\cdot), \mathfrak{E}(\cdot,\cdot),\mathfrak{S}(\cdot,\cdot) $ are the scalars. By allowing $m,n$ equal
to zero, we can simply express $w$ as
$$
\begin{aligned}
  w &=  \sum_{  ag_{\alpha}^m\in F(\Sigma)  }\mathfrak{E}(a, m) \lambda(a)\lambda(g_{\alpha})^m  \\
  w &=  \sum_{   g_{\alpha}^nb \in F(\Sigma)  }\mathfrak{S}(b, n)
   \lambda(g_{\alpha})^n\lambda(b)
\end{aligned}
$$
\noindent  Now we define the norms $||\ \cdot \ ||_{( \alpha, E)}$
and $||\ \cdot \ ||_{(  \alpha,  S)}$ as
$$\begin{aligned}
 & ||\ w \ ||_{( \alpha,  E)}    = \sqrt{  \sum_{g\in
 ES}|w(g)|^2
 } = \sqrt{ \sum_{ ag_{\alpha}^{m} \in ES}|\mathfrak{E}(a,m)|^2
 }\\
 & ||\ w \ ||_{(  \alpha,    S)}    = \sqrt{  \sum_{h\in SS}|w(h)|^2
 } = \sqrt{ \sum_{  g_{\alpha }^{ n}b  \in SS }
 | \mathfrak{S}(b,n)|^2
 }\end{aligned}
$$
 It is easy to see that, for every $w$ in $L(F(\Sigma))$, we have
\begin{equation}
\sum_{\alpha\in \Sigma} ||\ w \ ||_{( \alpha,  E)}^2\le
||w||_2^2,\qquad \sum_{\alpha\in \Sigma} ||\ w \ ||_{( \alpha,
S)}^2\le ||w||_2^2.\tag {\# }\end{equation}

\vspace{0.5cm}

\subsection{Another Estimation}

Following the notations from preceding subsection, we let $\Sigma$
be the index set; $\lambda$ be the left regular representation of
$F(\Sigma)$ and $L(F(\Sigma))$ be the free group factor with the
standard generators $\{\lambda(g_{\alpha})\}_{\alpha\in \Sigma}$.
For each $\alpha \in \Sigma$, let $y_{\alpha}$ be any self-adjoint
element in the von Neumann subalgebra generated by
$\lambda(g_{\alpha})$ such that $\tau(y_{\alpha})=0$.

 Fix some  index $\alpha  $    from $\Sigma$. Let $w_1, w_2$ be
two elements in $L(F(\Sigma))$. Then
$$
\begin{aligned}
  w_1 &=  \sum_{ ag_{\alpha}^m\in F(\Sigma)  }\mathfrak{E}_{w_1}(a, m) \lambda(a)\lambda(g_{\alpha})^m  \\
  w_2 &=  \sum_{ g_{\alpha}^nb \in F(\Sigma)  }\mathfrak{S}_{w_2}(b, n)
   \lambda(g_{\alpha})^n\lambda(b)
\end{aligned}
$$ where $\mathfrak{E}_{w_1}(\cdot,\cdot),\mathfrak{S}_{w_2}(\cdot,\cdot) $ are the scalars. Hence
$$
\begin{aligned}
 \tau(w_1y_\alpha w_2) & =\tau \left ( \sum_{ a, b ,m,n   }\mathfrak{E}_{w_1}(a, m) \lambda(a)\lambda(g_{\alpha})^m
  y_{\alpha} \mathfrak{S}_{w_2}(b,
 n)\lambda (g_{\alpha})^n\lambda (b)\right )\\
 &=\sum_{ a, b ,m,n   }\mathfrak{E}_{w_1}(a, m)\mathfrak{S}_{w_2}(b,
 n) \tau(y_{\alpha}\lambda (g_{\alpha})^{m+n}) \lambda(g_{\alpha})^n\lambda(b)\\
 &= I_1+I_2+I_3
\end{aligned}
$$
where
$$
 \begin{aligned}
  I_1& =\sum_{ a,  m=0 ,n\ne 0   }\mathfrak{E}_{w_1}(a, 0) \mathfrak{S}_{w_2}(a^{-1},
 n) \tau(y_{\alpha}\lambda(g_{\alpha})^{ n})  \\
   I_2& =\sum_{ a ,m\ne 0 ,n= 0   }\mathfrak{E}_{w_1}(a, m) \mathfrak{S}_{w_2}(a^{-1},
 0) \tau(y_{\alpha}\lambda(g_{\alpha})^{ m} ) \\
  I_3& =\sum_{ a  ,m\ne 0 ,n\ne 0   }\mathfrak{E}_{w_1}(a, 0) \mathfrak{S}_{w_2}(a^{-1},
 n) \tau(y_{\alpha}\lambda(g_{\alpha})^{ m+n}) \\
 \end{aligned}
$$

\begin{lemma} Let $y_{\alpha}$ be a   self-adjoint
element in the von Neumann subalgebra generated by Haar unitary
$\lambda(g_{\alpha})$ such that $\tau(y_{\alpha})=0$. We have the
following inequalities.
 $$
   \begin{aligned}
       |I_1| & \le ||w_1||_2 ||w_2||_{\alpha, S} ||y_{\alpha}||_2\\
       |I_2|  &\le ||w_1||_{\alpha, E} ||w_2||_2 ||y_{\alpha}||_2\\
        |I_3|
 & \le    ||w_1||_{\alpha, E}||w_2||_{\alpha,S}  \left ( \sum_k |\tau(y_{\alpha}\lambda(g_{\alpha})^{k})|  \right)
   \end{aligned}
 $$
\end{lemma}

\noindent {\bf Proof: } From Cauchy Schwartz inequality, we know
that
$$
  \begin{aligned}
    |I_1|  &\le \sum_{n\ne 0 }\left (
    \sqrt{\sum_{a} |\mathfrak{E}_{w_1}(a,    0)|^2}
        \sqrt{\sum_{a} |\mathfrak{S}_{w_2}(a^{-1}, n)|^2}
        \right ) |\tau(y_{\alpha}\lambda(g_{\alpha})^{ n})|\\
 & \le  \sqrt{\sum_{a} |\mathfrak{E}_{w_1}(a,
    0)|^2}  \sqrt{\sum_{a, n\ne 0} |\mathfrak{S}_{w_2}(a^{-1},
 n)|^2} \sqrt{\sum_{  n\ne 0} | \tau(y_{\alpha}\lambda(g_{\alpha})^{ n})|^2}\\
  & =  \sqrt{\sum_{a\in ET} |\mathfrak{E}_{w_1}(a,
    0)|^2}  \sqrt{\sum_{ g_{\alpha}^na^{-1}\in SS} |\mathfrak{S}_{w_2}(a^{-1},
 n)|^2}||y_{\alpha}||_2\\
 &\le ||w_1||_2 ||w_2||_{\alpha, S} ||y_{\alpha}||_2
  \end{aligned}
$$
Similarly, we have
$$
  \begin{aligned}
    |I_2|  \le ||w_1||_{\alpha, E} ||w_2||_2 ||y_{\alpha}||_2
  \end{aligned}
$$
And
$$
  \begin{aligned}
    |I_3|  & =\sum_{ a  ,k\ne 0 ,n\ne 0   }\mathfrak{E}_{w_1}(a, k-n) \mathfrak{S}_{w_2}(a^{-1},
 n) \tau(y_{\alpha}\lambda(g_{\alpha})^{ k})  \\
 &\le \sum_{k\ne 0 }\left (
    \sqrt{\sum_{ag_{\alpha}^n\in ES } |\mathfrak{E}_{w_1}(a,    k-n)|^2}
        \sqrt{\sum_{ g_{\alpha}^na^{-1} \in SS } |\mathfrak{S}_{w_2}(a^{-1}, n)|^2}
        \right ) |\tau(y_{\alpha}\lambda (g_{\alpha})^{ k})|\\
 & \le    ||w_1||_{\alpha, E}||w_2||_{\alpha,S}  \left ( \sum_k |\tau(y_{\alpha}\lambda(g_{\alpha})^{ k})|  \right)
  \end{aligned}
$$
Q.E.D.

\subsection{Proof of  Proposition 4.1}

\noindent {\bf Proof of Proposition 4.1: } Note that $B$ can be
written as $\sum_{1\le i,j \le n} b_{ij} \otimes e_{ij}$ with
$b_{ij}$ in $\N$ and
$$ ||E_n(B)||_2^2=||E_n(\sum_{1\le i,j\le n}b_{ij}\otimes e_{ij})||_2^2=\frac 1 n \sum_{1\le i,j\le n}|\tau_{\mathcal N}(b_{ij})|^2.  $$
Since the semicircular element $B$ and $\mathcal D_n$ are free with
respect to the trace $\tau_{\mathcal M}$, by Theorem 3.5, there are
unitary elements $u_1,\ldots, u_n$ in $\mathcal N$  such that
$UBU^*=\sum_{1\le i,j\le n}a_{ij}\otimes e_{ij}$ is a Voiculescu's
 semicircular matrix, where $U= \sum_{i} u_i\otimes e_{ii}$.
Therefore,
$$
B = \sum_{1\le i,j\le n} b_{ij}\otimes e_{ij}=U^*UBU^*U=\sum_{1\le
i,j\le n}(u_ia_{ij}u^*_j)\otimes e_{ij}, \text{ \ \ so \ \
$b_{ij}=u_ia_{ij}u_j^*$.  }
$$ where $\tau(a_{ij})=0$ and $\tau(a_{ij}a_{ij}^*)=  1 / n$.  We have
\begin{align}
 \frac 1 n \sum_{1\le i,j\le n}|\tau_{\mathcal N}(b_{ij})|^2&=
     \frac 1 n \sum_{1\le i \le n}|\tau_{\mathcal N}(b_{in})|^2 \\
     & \ +\frac 1 n \sum_{1\le  j\le n-1}|\tau_{\mathcal N}(b_{nj})|^2\\
     &\  + \frac 1 n \sum_{1\le i,j\le n-1}|\tau_{\mathcal N}(b_{ij})|^2
\end{align}
From Cauchy-Schwartz inequality, we have
$$\begin{aligned}
 (6) =\frac 1 n \sum_{1\le i \le n}|\tau_{\mathcal N}(b_{in})|^2 &=
     \frac 1 n \sum_{1\le i \le n}|\tau_{\mathcal
     N}(u_ia_{in}u_n)|^2\\
     &\le \frac 1 n \sum_{1\le i \le n}|\tau_{\mathcal
     N}(u_ia_{in}a_{in}^*u_i^*)\tau_{\mathcal N}(u_n^*u_n)|\\&= \frac 1 n \sum_{1\le i \le
     n} \frac 1 n  =\frac 1 n
\end{aligned}$$
Similarly, we have
$$ (7) = \frac 1 n \sum_{1\le  j\le n-1}|\tau_{\mathcal N}(b_{nj})|^2\le  \frac 1 n $$
Now we estimate (8),
\begin{align}
 (8) = \frac 1 n \sum_{1\le i,j\le n-1}|\tau_{\mathcal N}(b_{ij})|^2 &= \frac 1 n \sum_{1\le i,j\le n-1}|\tau_{\mathcal
  N}(u_ia_{ij}u_j^*)|^2 \notag \\
  &=\frac 1 n \sum_{1\le i \le n-1}|\tau_{\mathcal
  N}(u_ia_{ii}u_i^*)|^2 + \frac 2 n \sum_{1\le i<j\le n-1}|\tau_{\mathcal
  N}(u_ia_{ij}u_j^*)|^2 \notag\\
  &= 0 + \frac 2 n \sum_{1\le i<j\le n-1}|\tau_{\mathcal
  N}(u_ia_{ij}u_j^*)|^2=\frac 2 n \sum_{1\le i<j\le n-1}|\tau_{\mathcal
  N}(u_ia_{ij}u_j^*)|^2\notag\\
  &\le \frac 4 n \sum_{1\le i<j\le n-1} |\tau_{\mathcal
  N}(u_i x_{ij}  u_j^*)|^2 \\
  &+ \frac 4 n \sum_{1\le i<j\le n-1} |\tau_{\mathcal
  N}(u_i y_{ij}  u_j^*)|^2\\
  &\quad  \text {(where $x_{ij}, y_{ij}$ are the real and imaginal parts of $a_{ij}$.)
  }\notag
\end{align}
\noindent  Note that, since $UBU^*=\sum_{1\le i,j\le n}a_{ij}\otimes
e_{ij}$ is a Voiculescu's semicircular matrix,    $\{x_{ij},
y_{ij}\}_{1\le i<j\le n-1}$ is a free
 family of   semicircular
elements of $(0, r)$ with $r=\frac 1 {\sqrt {2n}} $.

  Let $\Sigma$ be the index set $\{\langle i, j \rangle\}_{1\le i< j\le n-1}$. Let $\mathcal N_1 $ be the von Neumann
subalgebra generated by $\{x_{ij}\}_{ 1\le i<j\le n-1 }$ in
$\mathcal N$, which is $*-$isomorphic to the free group factor $L(F(
\Sigma))$. Note that $\N_1$ is a subfactor of $\N$. By [10], there
exists a family of vectors $\{\xi_s\}_{s\in I} \subset L^2(\mathcal
N, \tau_{\mathcal N})$ such that, for each $u_i$,
$$
 u_i =\sum_{s}\xi_sE_{\mathcal N_1}(\xi_s^*u_i) = \sum_{s}\xi_sw(i,s),
$$ where $w(i,s) = E_{\mathcal N_1}(\xi_s^*u_i) \in \mathcal N_1$,
 and $E_{\mathcal N_1}(\xi_s^*\xi_t)=\delta_{st}f_s$, with $f_s$ projections in $\mathcal N_1$.
 And
it is easy to see that \begin{equation} \sum_{s}||f_sw(i,s)||_2^2
=||u_i||_2^2.\tag {\#\#}\end{equation}

\noindent Therefore,
$$
  \begin{aligned}
      \tau_{\mathcal N}(u_ix_{ij}u_j^*) &= \sum_{s_1,s_2}\tau_{\mathcal
      N}(\xi_{s_1}w(i,s_1)x_{ij}w(j,s_2)^*\xi_{s_2}^*)\\
      &= \sum_{s_1,s_2}\tau_{\mathcal
      N}(E_{\mathcal
      N_1}(\xi_{s_2}^*\xi_{s_1})(w(i,s_1)x_{ij}w(j,s_2)^*))\\
      &=\sum_{s_1}\tau_{\mathcal
      N} ((f_{s_1}w(i,s_1))x_{ij}(f_{s_1}w(j,s_1))^*))
  \end{aligned}
$$
For every  $\alpha= \langle i, j \rangle$, $w_1=f_{s_1}w(i,s_1)$,
$w_2= (f_{s_1}w(j,s_1))^*$, $x=x_{ij}$ and $y_{\alpha}= \tilde  x $,
we applied  Lemma 4.1 and 4.3, and obtained the following.
$$
\begin{aligned}
  |\tau_{\mathcal
      N} ((f_{s_1}w(i,s_1))\tilde x_{ij}(f_{s_1}w(j,s_1))^*))| &\le \left . \|f_{s_1}w(i,s_1)\|_2 ||f_{s_1}w(j,s_1)||_{(\langle i, j \rangle,S)}\cdot
      r\right .\\
&  +  \|f_{s_1}w(i,s_1)\|_{(\langle i, j \rangle,E)} ||f_{s_1}w(j,s_1)||_2\cdot r \\
 & + \left . ||f_{s_1}w(i,s_1)||_{(\langle i, j \rangle,E)} ||f_{s_1}w(j,s_1)||_{(\langle i,j\rangle,S)} \cdot \sqrt 5 r^{\frac 1
 2}\right. \\
\end{aligned}
$$
It follows that
$$
  \begin{aligned}
|\tau(u_i\tilde x_{ij}
   u_j^*)|& =|\sum_{s_1}\tau_{\mathcal
      N} ((f_{s_1}w(i,s_1))x_{ij}(f_{s_1}w(j,s_1))^*))|\\& \le \sum_{s_1}\left (\|f_{s_1}w(i,s_1)\|_2 ||f_{s_1}w(j,s_1)||_{(\langle i,j\rangle,S)}\cdot
      r\right .\\
&   +  \|f_{s_1}w(i,s_1)\|_{(\langle i, j \rangle,E)} ||f_{s_1}w(j,s_1)||_2\cdot r \\
 &   + \left . ||f_{s_1}w(i,s_1)||_{(\langle i, j \rangle,E)} ||f_{s_1}w(j,s_1)||_{(\langle i, j \rangle,S)} \cdot \sqrt 5 r^{\frac 1
 2}\right ) \\
  \end{aligned}
$$
By Lemma 4.2, we obtain
$$
  \begin{aligned}
    |\tau_{\mathcal N}(u_ix_{ij}u_j^*) |^2 & \le|\tau(u_i\tilde
      x_{ij}
   u_j^*)|^2 +6r^{\frac 5 2}.\\
       & \le \left |  \sum_{s_1}(\|f_{s_1}w(i,s_1)\|_2 ||f_{s_1}w(j,s_1)||_{(\langle i, j \rangle,E)}       \cdot
      r)\right . \\
&\quad \left . + \sum_{s_1}(\|f_{s_1}w(i,s_1)\|_{(\langle i, j \rangle,E)} ||f_{s_1}w(j,s_1)||_2\cdot r)\right . \\
 &\quad + \left . \sum_{s_1}(||f_{s_1}w(i,s_1)||_{(\langle i, j \rangle,E)} ||f_{s_1}w(j,s_1)||_{(\langle i, j \rangle,S)} \cdot \sqrt 5 r^{\frac 1
 2})\right |^2+ 6r^{5/2}\\
 &\le 3 \left |  \sum_{s_1}(\|f_{s_1}w(i,s_1)\|_2 ||f_{s_1}w(j,s_1)||_{(\langle i, j \rangle,S)} \cdot
      r)\right |^2 \\ 
      &\quad
 + 3 \left | \sum_{s_1}( ||f_{s_1}w(i,s_1)||_{(\langle i, j \rangle,E)}||f_{s_1}w(j,s_1)||_2\cdot r)\right |^2 \\
 &\quad + 3\left | \sum_{s_1}(||f_{s_1}w(i,s_1)||_{(\langle i, j \rangle,E)} ||f_{s_1}w(j,s_1)||_{(\langle i, j \rangle,S)} \cdot \sqrt 5 r^{\frac 1
 2}) \right |^2\\
 &\quad +   6r^{5/2}\\
      \end{aligned}
 $$  $$ \begin{aligned}
 &\le 3r^2 (\sum_{s_1}  ||f_{s_1}w(j,s_1)||_{(\langle i, j \rangle,S)}^2)(\sum_{s_1}  ||f_{s_1}w(i,s_1)||^2)\\&\quad +
 3r^2(\sum_{s_1}  ||f_{s_1}w(i,s_1)||_{(\langle i, j \rangle,E)}^2)(\sum_{s_1} ||f_{s_1}w(j,s_1)||^2)\\
 &\quad + 15r (\sum_{s_1} ||f_{s_1}w(j,s_1)||_{(\langle i.j\rangle , S}^2)(\sum_{s_1}
 ||f_{s_1}w(i,s_1)||_{(\langle i, j \rangle,E)} ^2)\\
 &\quad + 6r^{5/2}\\
 &\le 3r^2  \sum_{s_1}( ||f_{s_1}w(j,s_1)||_{{(\langle i, j \rangle,S)} }^2+
 ||f_{s_1}w(i,s_1)||_{(\langle i, j \rangle,E)} ^2)  \qquad \qquad \text {(because of
 (\#\#))}\\
 &\quad + 15r (\sum_{s_1}  ||f_{s_1}w(j,s_1)||_{(\langle i, j \rangle,S)} ^2)(\sum_{s_1}
 ||f_{s_1}w(i,s_1)||_{(\langle i, j \rangle,E)} ^2)\\
 &\quad +  6r^{5/2}\\
  \end{aligned}
$$
Therefore
$$\begin{aligned}
 (9) &\le \frac 4 n \sum_{1\le i<j\le n-1} |\tau_{\mathcal
  N}(u_i x_{ij}  u_j^*)|^2\\
    &\le \frac 4 n \sum_{1\le i<j\le n-1} \left ( 3r^2  \sum_{s_1} \left (||f_{s_1}w(j,s_1)||_{{(\langle i, j \rangle,S)} }^2+
 ||f_{s_1}w(i,s_1)||_{{(\langle i, j \rangle,E)} }^2\right )\right .\\
 &\left . \quad + 15r (\sum_{s_1}  ||f_{s_1}w(j,s_1)||_{{(\langle i, j \rangle,S)}}^2)(\sum_{s_1} ||f_{s_1}w(i,s_1)||_{{(\langle i, j \rangle,E)} }^2)
  +   6r^{5/2} \right)\\
 &\le \frac {12r^2}{n} \sum_{1\le j\le n-1} \sum_{s_1} \sum_{i}
 ||f_{s_1}w(j,s_1)||_{{(\langle i, j \rangle,S)} }^2     +  \frac {12r^2}{n} \sum_{1\le i\le n-1}
 \sum_{s_1}\sum_{j}
 ||f_{s_1}w(i,s_1)||_{{(\langle i, j \rangle,E)} }^2   \\
 &\quad  + \frac{60 r}{n} \sum_{1\le i\le n-1} \sum_{s_1}  \sum_{j}
 ||f_{s_1}w(i,s_1)||_{{(\langle i, j \rangle,E)} }^2   + \frac {24r^{5/2}}{n} \sum_{1\le i<j\le n-1} 1\\
 &\le  \frac {12r^2}{n} \sum_{1\le j\le n-1} 1+    \frac {12r^2}{n} \sum_{1\le i\le n-1} 1 +
 \frac{60 r}{n} \sum_{1\le i\le n-1} 1  + \frac {24r^{5/2}}
 {n} \sum_{1\le i<j\le n-1} 1\\
 &\quad\quad\quad\quad\quad\quad\quad\quad\qquad\qquad\qquad\qquad\quad \text {(because of }
 (\#) \text{ and } (\#\#) \text{)}
 \\
 &\le 24r^2+60r + 12 r^{5/2}n
 \end{aligned}
$$Similarly,
 we also have
 $$
   (10) \le  24r^2+60r + 12 r^{5/2}n.
 $$

\vspace{0.5cm}

\noindent Combining all of the above, we have
$$\|E_n(B)\|_2^2=\frac 1 n \sum_{1\le i,j\le n}|\tau_{\mathcal
N}(b_{ij})|^2\le \frac 2 n+ 48r^2+ 120r+ 24 r^{5/2}n\le   {49}
{\sqrt r}<\frac {49} {\sqrt [4] n} ,
$$ where $2r^2= \tau_{\mathcal N}(b_{ij}b_{ij}^*)= \frac 1 {n}$.
Q.E.D. \vspace{2cm}

\vspace{1cm}

\subsection{Definition of matricial distance of two elements}
\begin{Def}
Let $\mathcal M$ be a factor of type II$_1$ with the traical state
$\tau$. Let $\mathcal A$ be a diffused abelian von Neumann
subalgebra of
 $\mathcal M$.  Let $\mathcal S_n$ be the set
consisting of the type I$_n$ subfactors $\mathcal M_n$ of $\mathcal
M$ such that $\mathcal A\cap \mathcal M_n$ is a $n-$dimensional
subalgebra. Then, for any element $b$ in $\mathcal M$, the matricial
distance between $\mathcal A$ and $b$ is defined
as,$$\begin{aligned} \text{MatD}(\mathcal A,b)
=\liminf_{n\rightarrow \infty}\inf_{\mathcal M_n\in \mathcal S_n}\{
\|&b-E_{\mathcal M_n}(b)\|_2, \text{ where $E_{\mathcal M_n}$ is the
conditional
expectation}\\
 &\text { from $\mathcal M$ onto $\mathcal M_n$ }\}\end{aligned}$$
\end{Def}

\begin{Def}
  Let $a$ be a self-adjoint element in $\mathcal M$  such that $a$ generates a diffused abelian von
  Neumann subalgebra $\mathcal A$ of $\mathcal M$. Then, for any
  element $b$ in $\mathcal M$, the matricial distance between $a$
  and $b$ is defined as
  $$ MatD(a,b)= MatD(\mathcal A,b).   $$
\end{Def}

\begin{Def}
  Let $\mathcal A$, $\mathcal B$  be two diffused abelian von
  Neumann subalgebras of $\mathcal M$. Then, the matricial distance between $\mathcal A$
  and $\mathcal B$ is defined as
  $$ MatD(\mathcal A,\mathcal B)= \inf \{MatD(\mathcal A,b)/||b||_2\ | \  {b\in \mathcal B} \text { such that } E_{\A}(b)=0 \}.  $$
\end{Def}

\noindent From Proposition 4.1, we can easily have
\begin{Th}Let  $\M$ be a factor of type II$_1$.
Suppose that $a$ and $b$ are two free semicircular elements of
$(0,1)$ in   $\M$. Then $$ MatD (a,b) =1.$$
\end{Th}

\pb {\bf Question:} It is very interesting to consider the following
question: Suppose $\mathcal R$ is the hyperfinite II$_1$ factor and
$\mathcal A$ is a maximal abelian self-adjoint subalgebra of
$\mathcal R$. Does there exist some $x$ in $\mathcal R$ but not
contained in $\mathcal A$ such that $MatD(\mathcal A,x)=0? $

\pb {\bf Remark:} A positive answer to the preceding question,
combining with Theorem 4.1, will imply Popa's result that any one of
standard generator of free group factor $L(F(n))$ is not contained
in any hyperfinite II$_1$ subfactor of $L(F(n))$.

\pb {\bf Remark:} The further computation of matricial distance
between two free self-adjoint elements in a type II$_1$ factor will
be carried out in our forthcoming paper.

\newpage
\vspace{1cm} \noindent {\large\textbf{{Bibliography}}}

 \vspace{0.5cm}
  \begin{enumerate}





\item  L.\ Ge and S.\ Popa, ``On some decomposition properties for
factors of type II$_1$," Duke Math. J., 94 (1998), 79--101.

\item L. Ge ``On maximal injective subalgebras of factors," Adv. Math. 118 (1996), no. 1, 34--70

\item L. \ Ge and J. \ Shen, ``Applications of free entropy on
finite von Neumann algebras, III," \  GAFA, 12 (2002), no. 3,
546--566.


\item R.\ Kadison, ``Problems on von Neumann algebras," Paper
given at the Conference on Operator Algebras and Their Applications,
Louisiana State University, Baton Rouge, La., 1967

\item R.\ Kadison and J.\ Ringrose, ``Fundamentals of the Operator
Algebras,''  vols. I and II, Academic Press, Orlando, 1983 and 1986.

\item A. Nica; D. Shlyakhtenko; R. Speicher  ``$R$-cyclic families of matrices in free probability."
 J. Funct. Anal. 188 (2002), no. 1, 227--271.

\item A. Nica; R. Speicher  On the multiplication of free
$N$-tuples of noncommutative random variables. Amer. J. Math. 118
(1996), no. 4, 799--837.

\item S. Popa, ``On a problem of R. V. Kadison on maximal abelian
$*$-subalgebras in factors," Invent. Math. 65 (1981/82), no. 2,
269--281.

\item S. Popa, ``Maximal injective subalgebras in factors associated with free
groups," Adv. in Math. 50 (1983), no. 1, 27--48.

\item S. Popa, ``Classification of subfactors of finite depth of the
hyperfinite type ${\rm III}\sb 1$ factor," C. R. Acad. Sci. Paris
S¨¦r. I Math. 318 (1994), no. 11, 1003--1008.



\item D. Voiculescu, ``Circular and semicircular systems and free product factors," Operator algebras,
unitary representations, enveloping algebras, and invariant theory
(Paris, 1989), 45--60, Progr. Math., 92, Birkhauser Boston, Boston,
MA, 1990.

\item   D. Voiculescu, ``The analogues of entropy and of Fisher's
information measure in free probability theory II," Invent.\ Math.,
{  118} (1994), 411-440.

\item D. Voiculescu, ``The analogues of entropy and of Fisher's
information measure in free probability theory III: The absence of
Cartan subalgebras," Geom.\ Funct.\ Anal.\ { 6} (1996) 172--199.

\item D. Voiculescu, K. Dykema and A. Nica,  { ``Free Random
Variables,''} CRM Monograph Series, vol. 1, AMS, Providence, R.I.,
1992.


\end{enumerate}

\end{document}